\documentclass[12pt]{article}
 \usepackage{amsfonts}

 \newcommand{\oper}[2]{\newcommand{#1}{\mathop{\mathrm{#2}}\nolimits} }
 \oper{\tr}{tr}
 \oper{\adj}{adj}
 \oper{\Div}{div}
 \oper{\ad}{ad}
 \oper{\Ad}{Ad}
 \oper{\End}{End}
  \oper{\Hom}{Hom}
  \oper{\Aut}{Aut}
  \oper{\SO}{SO}
  \oper{\SP}{Sp}
  \oper{\SU}{SU}
  \oper{\GL}{GL}
  \oper{\T}{T}
  \oper{\U}{U}
  \oper{\id}{I}
  \oper{\ext}{Ext}
  \oper{\rank}{rank}
  \oper{\diag}{Diag}

  \def\bC{{\bf C}}
  
  \def\bR{{\bf R}}
  \def\bH{{\bf H}}
  \def\bZ{{\bf Z}}

  \def\cp{{\bf C}{\bf P}}

  \def\cd{{\cal D}}
  \def\ce{{\cal E}}

  \def\co{{\cal O}}

  \def\cv{{\cal V}}

  \def\dbar{\overline\partial}
  \def\wbar{\overline w}
  \def\zbar{\overline z}

  \def\mubar{\overline\mu}

  \def\osigma{\overline\sigma}
  \def\oomega{\overline\omega}
  \def\Oomega{\overline\Omega}

  \def\Up{\Upsilon}

  \newcommand{\lie}[1]{\mathfrak{#1}}
  \newcommand{\pdd}[1]{\frac{\partial}{\partial{#1}}}

  \def\hei{\check{H}}

  \newtheorem{corollary}{Corollary}
  \newtheorem{proposition}{Proposition}
  \newtheorem{theorem}{Theorem}
  
  \newtheorem{lemma}{Lemma}
  \newtheorem{remark}{Remark}

  \newcommand{\bproof}{\noindent{\it Proof: }}
  \newcommand{\eproof}{\  q.~e.~d. \vspace{0.2in}}

\begin{document}

  \title{Deformations of Hypercomplex Structures\\
  associated to Heisenberg Groups}
  \author{Gueo Grantcharov
  \thanks{
  Permanent Address: Department of Mathematics, Florida International University, Miami, FL 33199,
  Partially supported by NSF DMS-0333172}
  \and
  Henrik Pedersen
  \thanks{Address: Department of Mathematics and Computer Science, University of Southern Denmark, Campusvej 55, Odense M, DK-5230,
       Denmark. E-mail: henrik@adm.sdu.dk. Partially supported by the European
       contract HPRN-CT-2000-00101.}
  \and
  Yat Sun Poon
  \thanks{
  Address: Department of Mathematics, University of California at Riverside,
  Riverside, CA 92521, U.S.A.. E-mail: ypoon@math.ucr.edu. Partially supported by NSF DMS-0204002.}
  }
  %\date{hei/020918.tex, Geo's addition of moduli}
  \maketitle

  \noindent{{\bf Abstract}: Let $X$ be a compact quotient of the product
  of the real Heisenberg group $H_{4m+1}$ of dimension $4m+1$ and the 3-dimensional
  real Euclidean space $\bR^3$. A left invariant hypercomplex structure
  on $H_{4m+1}\times \bR^3$ descends onto the compact quotient $X$.
  The space $X$ is a hyperholomorphic fibration of 4-tori over
  a $4m$-torus.  We calculate
  the parameter space and obstructions to  deformations of this hypercomplex
  structure on $X$. Using our calculations we show that
all small deformations generate invariant hypercomplex structures
on $X$ but not all of them arise from deformations of the lattice.
This is in contrast to the deformations on the $4m$-torus.

  \

  \noindent{Keywords: Heisenberg group, nilpotent groups, Kodaira surfaces,
  Kodaira manifolds, hypercomplex, deformation.}

  \

  \noindent{AMS Subject Classification: Primary 32G07. Secondary 53C15,
  53C56, 32L25, 57S25.}

  \section{Introduction}

  Heisenberg groups play a fundamental role in many branches of mathematics.
  One of its appearances is in the construction of Kodaira surfaces.
  Among other features, these surfaces can be realized as an elliptic fibration
   over elliptic curves. This is a realization of the quotient
   map from the Heisenberg group with respect to its center. This construction
  is extended to higher dimensions in \cite{GPP} \cite{MPPS}. On the other hand,
  Heisenberg groups are used by several authors to construct
  hypercomplex structures \cite{Joyce} \cite{Dotti} \cite{Laura}.
  This direction greatly enriches the source of hypercomplex
  manifolds as the past constructions of compact examples
  are often limited to homogeneous spaces with semi-simple Lie groups
  \cite{Joyce2} \cite{PP2}.
  The construction and the deformations of these hypercomplex
  structures are the topics of this article.

  Let  $H_{2n+1}$ be the $2n+1$-dimensional real Heisenberg group.
  It   is the  extension of the Abelian additive group
  $\bR^{2n}$ by a one-dimensional center.
  The space $\bR^{2n}$ is a complex space as it can be real linearly identified to
  the complex vector space $\bC^n$.
  The product $H_{2n+1}\times \bR^1$ admits a left-invariant complex structures
  such that the natural projection $\phi$  onto $\bR^{2n}$ is holomorphic.
  Taking compact quotients, we obtain a complex structure on a generalization
  of \it Kodaira manifolds \rm $X$ \cite{GPP}. This construction can be
  extended further to a construction of hypercomplex structure.

  The quotient of
  the $4m+1$-dimensional real Heisenberg group
   $H_{4m+1}$ with its one-dimensional center is the Abelian additive group
  $\bR^{4m}$.
  The space $\bR^{4m}$ is a hypercomplex space when it is identified to
  the module of quaternions $\bH^m$.
  The product $H_{4m+1}\times \bR^3$ admits a left-invariant hypercomplex structures
  such that the natural projection $\phi$  onto $\bR^{4m}$ is hyper-holomorphic.
  Taking compact quotients $X$, we obtain a hypercomplex manifold fibered over
  the torus $T^{4m}$ with its standard hypercomplex structures.
  The fiber is a four-dimensional torus $T^4$ obtained as a compact quotient
  of the product of the group $H_{4m+1}\times \bR^3$. In this paper,
  we first calculate the parameter spaces and obstructions to
  deformation of the hypercomplex structure on the compact quotient $X$ of
  $H_{4m+1}\times \bR^3$. Our computation is based on twistor theory, deformation
  theory of maps and an understanding of the deformation of hypercomplex structures
  on the torus. Then we compare this parameter space with the
  space of the invariant hypercomplex structures and the
  deformation space arising from the deformations of the lattice in $H_{4m+1}\times \bR^3$.

  In Section 1, we explain the relations among twistor theory,
   deformation theory of holomorphic maps and deformation theory
  of hypercomplex manifolds.
  Through these relations, we
  compute the parameter space for hypercomplex deformations in
  Section 2. The next theorem is the result of enumerating the dimension
  of the parameter space in a long exact sequence of cohomology.
  \begin{theorem}\label{theorem1}
   The real dimension of the virtual parameter space of deformations of
  hypercomplex structures on the $(4m+4)$-dimensional
  manifold $X$ is equal to $6m^2+11m+12$.
  \end{theorem}
  In the parameter space, there is a twelve-dimensional subspace contributed
  by the deformation of the 4-torus in the fiber of the projection $\phi$.
  Hypercomplex deformations of the base of this projection contribute
  to a subspace of dimension $3(2m^2+m)$.
  In Lemma \ref{torus moduli}, we establish that
  the dimension of the parameter space of hypercomplex structures on
  a torus of dimension $4m$ is naturally identified to
  $12m^2$ and that only some of the hypercomplex deformation on the base
  torus comes from a hypercomplex
  deformation on the Kodaira manifold.
  Since the space of obstructions to deformation does not vanish, we study the integrability of deformation parameters by constructing
  convergent power series in Section 3.
  \begin{theorem}\label{theorem2} Every point in the virtual parameter space is an infinitesimal
  deformation of an integrable deformation.
  \end{theorem}
  In this construction, we do not control the power series enough to claim a priori that the
  deformation must be hypercomplex. At this point we produce a deformation of quaternionic
  structures only. We conclude our work showing that the deformed twistor
  spaces have a holomorphic projection onto $\cp^1$. Through the
  twistor correspondence, we complete a proof of
  the following theorem.
  \begin{theorem}\label{theorem3} Every quaternionic deformation of the hypercomplex structure
  on $X$ is a hypercomplex deformation.
  \end{theorem}
  This observation raises the issue of enumerating the number of parameters for quaternionic
  deformations of the hypercomplex manifold $X$. This is achieved through a coboundary
  map computation.
  \begin{theorem}\label{theorem4} The real dimension of the parameter space of deformations of quaternionic structures
 on the $(4m+4)$-dimensional manifold $X$ is equal to $6m^2+
  11m+9$.
  \end{theorem}

  During the course of our computation, we reveal a considerable amount of
  information about the deformation of hypercomplex structures on the torus $T^{4m}$,
  because its twistor space $Z$ is the base space of an elliptic fibration from
  the twistor space $W$ of the hypercomplex structure $X$. In section 4.5
 we reveal how our computation may be used to identify the moduli space
  of hypercomplex structures on a torus.

Finally in Section 5, we calculate the parameter space of the
invariant hypercomplex structures on $X$ and space of deformations
arising from the deformations of the lattice. Although in Section
3 we do not control the power series completely, the information
we obtain there is enough to prove the first part of the
following:

  \begin{theorem}\label{theorem5}
  Any small deformation of the hypercomplex structure on $X$ consists of invariant structures. However
  there are deformations which
  do not arise from a deformation of the lattice $\Gamma$.
  \end{theorem}

The second part of Theorem 5 follows by direct dimension count.
The proof is at the end of Section 5.

  \section{Basic Constructions}

  \subsection{The Heisenberg group and the Heisenberg algebra}\label{algebra}
  The real \it Heisenberg group \rm of dimensional 4m+1 is the Lie group $H_{4m+1}$
  whose underlying manifold is $\bR^{4m}\times \bR$ with
  coordinates $(x,y,z):=(x_1, \dots, x_{2m}, y_1, \dots, y_{2m}, z)$ and whose group
  law is given by
  \begin{equation}\label{group law}
  (x, y, z)*(x', y', z')=
  (x+x', y+y', z+z'-2\sum_{j=1}^{2m}(x_jy_j'-y_jx_j')).
  \end{equation}
  The left translations of $\{{\frac{\partial}{\partial x_j}}_{|{\mathbf 0}},
   {\frac{\partial}{\partial y_j}}_{|{\mathbf 0}},
   {\frac{\partial}{\partial z}}_{|{\mathbf 0}}\}$
  are the following vector fields.
  \begin{equation}\label{fields}
  X_j=\frac{\partial}{\partial x_j}+2y_j\frac{\partial}{\partial z},
  \hspace{.2in}
  Y_j=\frac{\partial}{\partial y_j}-2x_j\frac{\partial}{\partial z},
  \hspace{.2in}
  Z=\frac{\partial}{\partial z}.
  \end{equation}
  These vectors form a basis for
  the \it Heisenberg algebra \rm $\lie{h}_{4m+1}$
  of the Heisenberg group $H_{4m+1}$.
  The commutation relations are as follows:
  \begin{equation}
  [Y_j, X_k]=4\delta_{jk} Z,
  \hspace{.2in}
  [X_j, X_k]=[Y_j, Y_k]=[X_j, Z]=[Y_j, Z]=0.
  \end{equation}
  The subspace $\lie{c}$ spanned by $Z$ is the center of
  the Heisenberg algebra. The quotient space of the Heisenberg algebra
  with respect to the center is the $4m$-dimensional Abelian algebra
  $\lie{t}_{4m}$. Therefore, we have the exact sequence
  \begin{equation}\label{extension}
  0\to \lie{c} \stackrel{\iota}{\rightarrow}
    \lie{h}_{4m+1} \stackrel{\phi}{\rightarrow} \lie{t}_{4m} \to 0.
  \end{equation}
  On the level of Lie groups, we have a group homomorphism
  \begin{equation}
  \phi: H_{4m+1}=(\bR^{4m+1}, *) \rightarrow (\bR^{4m}, +)
  \end{equation}
  obtained as the quotient of the central subgroup $C=(\bR^1, +)$.
  Although it is obvious that
  \begin{equation}\label{lifting}
  \phi(X_j)=\frac{\partial}{\partial x_j},
  \qquad
  \phi(Y_j)=\frac{\partial}{\partial y_j},
  \qquad
  \phi(Z)=\frac{\partial}{\partial z},
  \end{equation}
  it will be important for our future computation that these
  identities give a way to lift vector fields from the Abelian
  group to the Heisenberg group.

  Let $\Gamma$ be the subgroup $(\bZ^{4m+1}, *)$ of the Heisenberg group.
  The intersection $\Gamma_0$ of $\Gamma$ with the central subgroup
  is isomorphic to the integer group $\bZ$.
  The quotient of $C$ by $\Gamma_0$ is the one-dimensional torus group $T^1$.
  The quotient of the additive Abelian group $\bR^{4m}$ by $\phi(\Gamma)$ is the
  4m-dimensional torus group $T^{4m}$.
  It is obvious that the homomorphism $\iota$ from
  the central subgroup $C$ into the Heisenberg group $H_{4m+1}$ intertwines
  $\Gamma_0$ and $\Gamma$, and the homomorphism $\phi$ from the
  Heisenberg group $H_{4m+1}$ to the Abelian group $\bR^{4m}$ intertwines
  the groups $\Gamma$ and $\phi(\Gamma)$. Therefore, the homomorphisms $\iota$ and
  $\phi$ descend to maps between compact quotients. Since $C$ is the central
  subgroup, its action commutes with the action of the
  lattice $\Gamma$. Therefore, the quotient group $T^1$ acts on the quotient
  space of left-cosets $\Gamma\backslash {H}_{4m+1}$.
  The orbits of this group action are precisely the
  fibers of the projection
  \begin{equation}
  \phi: \Gamma\backslash {H}_{4m+1}\rightarrow T^{4m}.
  \end{equation}          From now on,
  we denote the quotient space $\Gamma\backslash H_{4m+1}$ by
  ${\check{H}}_{4m+1}$ or $\hei$ if the dimension of the group is clear in
  a given context.

  \subsection{A Construction of Hypercomplex Structures}\label{defining}
  Three complex structures $I_1, I_2$ and $I_3$ on a smooth manifold form a
  hypercomplex structure if
  \begin{equation}\label{quaternion}
  I_1^2=I_2^2=I_3^2=-I_0,
  \quad
  \mbox{ and }
  \quad
  I_1I_2=I_3=-I_2I_1,
  \end{equation}
  where $I_0$ is the identity map.

  Let $\lie{t}_3$ be the
  3-dimensional Abelian algebra.
  The direct sum  $\lie{h}_{4m+1} \oplus \lie{t}_3$ is a 2-step
  nilpotent algebra whose center is four-dimensional.
  Fix a basis $\{ E_1, E_2, E_3\}$ for
  $\lie{t}_3$. Consider the endomorphisms $I_1$, $I_2$ and $I_3$
  of $\lie{h}_{4m+1} \oplus \lie{t}_3$ defined by  left multiplications of
  the quaternions $i$, $j$ and $k$ on the module of quaternions $\bH$, and
  the identifications
  \begin{eqnarray}\label{identification}
  x_{2a-1}X_{2a-1}+x_{2a}X_{2a}+y_{2a-1}Y_{2a-1}+y_{2a}Y_{2a}
  &\rightarrow &x_{2a-1}+x_{2a}i+y_{2a-1}j+y_{2a}k;
  \nonumber\\
  zZ+t_1E_1+t_2E_2+t_3E_3
  &\rightarrow &z+t_1i+t_2j+t_3k.
  \label{hyper}
  \end{eqnarray}
  In other words, for $1\leq a\leq m$,
  \begin{eqnarray*}
  I_1X_{2a-1} = X_{2a}, &
  I_1Y_{2a-1}=Y_{2a}, &
  I_1Z=E_1, \hspace{.2in} I_1E_2=E_3;\\
  I_2X_{2a-1} = Y_{2a-1}, &
  I_2X_{2a}=-Y_{2a}, &
  I_2Z=E_2,  \hspace{.2in} I_2E_1=-E_3;\\
  I_3X_{2a-1}=Y_{2a}, & I_3X_{2a}=Y_{2a-1}, &
  I_3Z=E_3, \hspace{.2in} I_3E_1=E_2.
  \end{eqnarray*}
  Through left translations, these endomorphisms define almost complex structures
  on the product  of the Heisenberg group and the three-dimensional
  additive group $H_{4m+1}\times \bR^3$.
  By construction, these almost complex
  structures satisfy the algebra (\ref{quaternion}).
  As $[I_aX, I_aY]=[X, Y]$  for any
  left-invariant vector fields $X$ and $Y$ and $1\leq a\leq 3$,
  these complex structures are integrable.
  It implies that $ \{ I_a : a =1, 2, 3 \} $
   is a left-invariant hypercomplex structure
  on the Lie group $H_{4m+1}\times \bR^3$. Let $\bZ^3$ be the integer subgroup
  of $\bR^3$.
  Then the quotient space $\hei\times \bR^3/\bZ^3$
   is the  compact  manifold $\hei\times T^3$. We denote this compact
  hypercomplex manifold by $X$.
  The natural projection $\phi$  from $X$ onto $T^{4m}$
  and the inclusion $\iota$ from $T^1\times T^3\cong T^4$ into
  $X$ are both hyper-holomorphic maps. The fibers of the projection
  $\phi$ are the orbits of the left-action of $T^1\times T^3$.

 For the hypercomplex structure we show:

\begin{theorem}\label{identification} The hypercomplex structure on $H_{4m+1}\times {\bf R}^3$ constructed above
is  equivalent to the standard one on ${\bf H}^{m+1}$.
\end{theorem}

{\it Proof:} From the definition of the vector fields
$X_i,Y_i,Z,E_i$ we see that the dual 1-forms are $dx_i,dy_i,\theta
= dz-2\Sigma(y_idx_i-x_idy_i),de_i$, where $e_1,e_2,e_3$ are
coordinates of ${\bf R}^3$. Then we have the same identities for
the action of the hypercomplex structure on the 1-forms as we had
for the vector fields, e.g. $I_1dx_1=dx_2...$ etc. With this in
mind we calculate:
$$
\begin{array}{lll}
I_1dz &=& I_1\theta+I_12\sum_i(y_idx_i-x_idy_)=\\
&=&de_1+2\sum_a(y_{2a-1}dx_{2a}-x_{2a-1}dy_{2a}-y_{2a}dx_{2a-1}+x_{2a}dy_{2a-1})\\
&=&d(e_1+2\sum_a(y_{2a-1}x_{2a}-x_{2a-1}y_{2a}))
\end{array}
$$
so by defining
$$f_1:=e_1+2\sum_a(y_{2a-1}x_{2a}-x_{2a-1}y_{2a})$$ we have
$I_1dz=df_1$. Similarly,
$$
\begin{array}{lll}
I_2dz &=& I_2\theta+I_22\sum_i(y_idx_i-x_idy_)=\\
&=&de_2+2\sum_a(y_{2a-1}dy_{2a-1}+x_{2a-1}dx_{2a-1}-y_{2a}dy_{2a}-x_{2a}dx_{2a}))\\
&=&d(e_2+\sum_a(y_{2a-1}^2+x_{2a-1}^2-y_{2a}^2-x_{2a}^2))
\end{array}
$$
so $I_2 dz=df_2$, for
$$f_2:=e_2+\sum_a(y_{2a-1}^2+x_{2a-1}^2-y_{2a}^2-x_{2a}^2)$$
Finally, by similar calculation we have $I_3 dz=df_3$, where
$$
f_3=e_3+2\sum_a(y_{2a-1}y_{2a}+x_{2a-1}x_{2a})
$$

Now from here we have that $x_i,y_i,z,f_i$ are quaternionic coordinates,
which are global on ${\bf R}^{4m+4}$, so identify our hypercomplex structure with
the standard one as claimed. $q.e.d.$

The above calculation could be done using the Obata connection and then one can
 show that $df_i$ are parallel 1-forms.
Here we outline the argument.
The Obata connection is given by
$$\nabla_XY=1/2[X,Y]+1/12\sum_{i,j,k}I_i([I_jX,I_kY]+[I_jY,I_kX])
+1/6\sum_i I_i([I_iX,Y]+[I_iY,X])$$ for hypercomplex manifold,
where $(i,j,k)$ is cyclic permutation of $(1,2,3)$
 . Then for abelian hypercomplex structure this reduces to
$$\nabla_X Y= 1/2[X,Y]+1/2\sum_i I_i[I_iX,Y]$$

From here one has that  $dx_i,dy_i$ are parallel and
% $\theta$ all $\nabla_X\theta(Y)$
%vanish except $\nabla_{X_1}\theta(Y_1)=\nabla_{X_2}\theta(Y_2)=2,
%\nabla_{Y_1}\theta(X_1)=\nabla_{Y_2}\theta(X_2)=-2$.
 $\nabla\theta = 1/2d\theta$. Now from $\nabla(x_idy_j)=dx_i\otimes dy_j$ one
easily checks that $dz$ and $df_i$ are parallel.
% Moreover from the fact that $I_j$ commute with $\nabla$, one calculates
%$\nabla de_i$ and finds $f_j$ above.

  \subsection{Twistor Theory}

  We identify points
  $\vec{a}=(a_1, a_2, a_3)$ in the unit 2-sphere $S^2$ to the complex structure
  $I_{\vec{a}}=a_1I_1+a_2I_2+a_3I_3$. Let $J_{\vec{a}}$ be the complex
  structure on $S^2$ defined by the stereographic projection
  \begin{equation}\label{stereographic}
  \mu\in \bC \mapsto {{\vec{a}}}
  =\frac{1}{1+|\mu |^{2}}
  (|\mu |^{2}-1, -i(\mu -{\overline{\mu }}), \mu +{\overline{\mu }}).
  \end{equation}
  This map takes $0$ to $-I_1$, $i$ to $I_2$ and $1$ to $I_3$.
  It sends the complex orientation of the complex
  plane to the outward normal orientation of the sphere.
  Therefore, the complex structure on the sphere at the unit vector $\vec{a}$ is
  defined by the cross product with $\vec{a}$.

  The smooth manifold $W=X\times S^2$
  is endowed with an almost complex structure
  $\cal{I}$ defined by ${\cal{I}}_{(x, {\vec{a}})}=I_{\vec{a}}\oplus J_{\vec{a}}$.
  By twistor theory, this is an integrable complex structure \cite{PP1}. Moreover,
  the projection $p$ from $W$ onto $S^2$ is a holomorphic projection
  such that the fiber $p^{-1}(\vec{a})$ is the complex manifold $(X, I_{\vec{a}})$.
    The holomorphic projection from $Z$ onto $\cp^1$ is also denoted by $p$.
  The map is also \it real \rm in the sense that there is an anti-holomorphic
  involution $\tau$ on the twistor space $W$ such that $p\circ\tau=\rho\circ p$ where
  $\rho$ is the anti-podal map on the 2-sphere.

  As explained in \cite{PP2}, deformations of hypercomplex structures
  are identified to deformations of the real map $p$.
  Deformations of this map $p$ are described by the cohomology spaces
  $H^k(W, \cd_W)$ where
   $\cd_W$ is the kernel of the differential $dp$
    \cite{Hor}. The real part of these spaces contain the deformation
  of the hypercomplex structures. Since there is also a correspondence
  between quaternionic structure and the complex structures on the twistor space
  \cite{Salamon}, the real part of the cohomology spaces $H^k(W, \Theta_W)$
  contains the deformation theory of the quaternionic structures.

  \section{Deformation Parameters of $X=\hei\times T^3$}

   The aim of this section is to compute the cohomology spaces
  $H^k(W, \cd_W)$ of the
  twistor space $W$ for the hypercomplex manifold $X$.
  We begin with some standard computation on the twistor space $Z$ over
  the torus.

  \subsection{Twistor Space of the Torus}
  The twistor space $Z$ is the quotient of the bundle
  $\lie{t}^{1,0}\otimes \co (1)$ on $\cp^1$  where
  $\lie{t}^{1,0}$ is the $(1,0)$-part of the complexification of
  the $4m$-dimensional Abelian algebra $\lie{t}$
  \cite[Example 13.64 and Example 13.66]{Besse}.

  Choose linear coordinates $(z_1^a, z_2^a, z_1, z_2), a=1, \dots, m,$
  for $\bC^{2m+2}$.
  They are related to real coordinates of $H_{4m+1}\times\bR^3$ by
  \begin{equation}
  z_1^a=x_{2a-1}+ix_{2a},
  \hspace{.2in}
  z_2^a=y_{2a-1}+iy_{2a},
  \hspace{.2in}
  z_1=z+it_1,
  \hspace{.2in}
  z_2=t_2+it_3.
  \end{equation}

  Let $[\lambda_1, \lambda_2]$ be the homogeneous coordinates on $\cp^1$.
  On $U_1=\{\lambda\in\cp^1: \lambda_1\neq 0\}$, define $\nu=\frac{\lambda_2}{\lambda_1}$.
  On $U_2=\{\lambda\in\cp^1: \lambda_2\neq 0\}$, define $\mu=\frac{\lambda_1}{\lambda_2}$.
  We use the same notation to denote $p^{-1}(U_1)$ and $p^{-1}(U_2)$ on both
  $Z$ and $W$.

  \begin{lemma}\label{direct}
   Let $\co$ be the structure sheaf of the twistor space $Z$.
  Let  $R^q p_*\co$ be the q-th direct image sheaf with respect to the projection $p$
  from $Z$ onto $\cp^1$. Then
  \begin{equation}
  R^q p_*\co=\wedge^q \left( \lie{t}^{*(0,1)}\otimes \co (1)\right).
  \end{equation}
  \end{lemma}
  \bproof
  Since the dimension of $H^q(p^{-1}(\lambda ), \co )$ is constant with respect to $\lambda$,
   the direct image sheaves $R^q p_*\co$ are locally free.
  As each fiber $p^{-1}(\lambda )$ is an Abelian variety,
  $
  H^q(p^{-1}(\lambda ), \co )=\wedge^q H^1(p^{-1}(\lambda ), \co),
  $
  and there is a vector bundle isomorphism $R^q p_*\co=\wedge^q R^1p_*\co$.

  On $U_2$,  for $1\leq a\leq m$,  the (0,1)-forms
  \begin{equation}\label{fourteen}
  {\overline\sigma}_1^a=\frac{\mubar d\zbar_1^a-d{z}_2^a}{1+|\mu |^2},
  \hspace{.2in}
  {\overline\sigma}_2^a=\frac{\mubar d\zbar_2^a+d{ z}_1^a}{1+| \mu |^2}
  \end{equation}
  are holomorphic because
  \[
  d{\overline\sigma}_1^a=\frac{1}{1+|\mu |^2}
      \left( d\mubar\wedge{\sigma}_2^a
        -\mubar d{\mu}\wedge \osigma_1^a
      \right),
  \hspace{.2in}
  d{\overline\sigma}_2^a=\frac{1}{1+|\mu |^2}
      \left(-\mubar d{\mu}\wedge{\overline\sigma}_2^a
       -d\mubar\wedge{\sigma}_1^a \right)
  \]
  are type (1,1)-forms. Since for every $\lambda$, $h^1(p^{-1}(\lambda ), \co)=2m$,
  these holomorphic forms determine a trivialization of the bundle
  $R^1p_*\co (U_2)$.
  Similarly, on $U_1$ the (0,1)-forms
  \begin{equation}
  {\overline\rho}_1^a=\frac{d\zbar_1^a-{\overline\nu} d{z}_2^a}{1+|\nu |^2},
  \hspace{.2in}
  {\overline\rho}_2^a=\frac{d\zbar_2^a+{\overline\nu} d{z}_1^a}{1+| \nu |^2}
  \end{equation}
  determine a holomorphic trivialization of the bundle $R^1p_*\co (U_1)$.
  As ${\overline\rho}_i^a=\mu {\overline\sigma}_i^a$,
   the bundle $R^1p_*\co$ is isomorphic to
  $\lie{t}^{*(0,1)}\otimes \co (1)$ as claimed.
  \eproof

  \begin{lemma}\label{sections}
  For $\ell\geq -1$,
  $H^{k}(Z, p^*\co (\ell ))=
  \lie{t}^{*(0,k)}\otimes S^{\ell+k}\bC^{2}$
  where $S^{j}\bC^2$ is the j-th symmetric tensor product of $\bC^2$.
  \end{lemma}
  \bproof
  The projection formula and the last lemma determine the isomorphism
  \[
  R^qp_*p^*\co (\ell ) = R^qp_*\co \otimes\co (\ell )
      =\co (\ell )\otimes \wedge^q
    \left( \lie{t}^{*(0,1)}\otimes \co (1)\right)
      =\co (\ell+q)\otimes \lie{t}^{*(0,q)}.
  \]
  Consider the Leray spectral sequence with
  $E_2^{p,q}=H^p(\cp^1, R^qp_*p^*\co (\ell )),$ and
  $E_{\infty}^{p, q}\Rightarrow H^{p+q}(Z, p^*\co (\ell )).$
  When $\ell \geq -1$,
  $E_2^{p, q}=0$ for all $p\geq 1$ and $q\geq 0$. Therefore,
  the spectral sequence degenerates at $E_2$, and
  \begin{eqnarray*}
  H^{k}(Z, p^*\co (\ell )) &=&\oplus_{p+q=k}E_2^{p, q}
       = E_2^{0,k}=H^0(\cp^1, R^kp_*p^*\co (\ell ))\\
     &=& H^0(\cp^1, \co (\ell+k)\otimes \lie{t}^{*(0,k)})
     = \lie{t}^{*(0,k)}\otimes S^{\ell+k}\bC^{2}.
  \end{eqnarray*}
  \eproof

  The cohomology spaces of the last lemma can be described explicitly.
  Define
  \begin{equation}\label{omegas}
  \Oomega_1^a=\frac
  {{\overline\lambda}_1d{\overline z}_1^a-{\overline\lambda}_2dz_2^a}
  {|\lambda_1|^2+|\lambda_2|^2},
  \hspace{.2in}
  \Oomega_2^a=\frac
  {{\overline\lambda}_1d{\overline z}_2^a+{\overline\lambda}_2dz_1^a}
  {|\lambda_1|^2+|\lambda_2|^2}.
  \end{equation}
  Then $\{\lambda_1\Oomega_1^a, \lambda_2\Oomega_1^a,
  \lambda_1\Oomega_2^a, \lambda_2\Oomega_2^a\}$ forms a basis for
   the space $H^1(Z, \co )$. More generally,
  the space $H^k(Z, p^*\co(\ell ))$, for $\ell \geq 0$ is spanned by
  the twisted k-forms
  \begin{equation}\label{basic}
    \left(\lambda_1^{\ell+k-l}\lambda_2^l\right)
    \Oomega_1^{a_1}\wedge\cdots\wedge\Oomega_1^{a_i}
        \wedge\Oomega_2^{b_1}\wedge\cdots\wedge\Oomega_2^{b_j}
  \end{equation}
  where $0\leq l\leq \ell+k$.
  The space $H^1(Z, \co)$ has an alternative description.
  For $k=0,1,2,3$ and
  over $p^{-1}({\vec{a}})$, define 1-forms
  \begin{equation}\label{twisted forms}
  {\overline\omega}_k^a=I_kdx_{2a-1}-iI_{\vec{a}}I_kdx_{2a-1}.
  \end{equation}
  These (0,1)-forms on the twistor space span the space $H^1(Z, \co)$
  because
  \begin{eqnarray}\label{forms}
  \oomega_0^a=\mu {\overline\sigma}_1^a+{\overline\sigma}_2^a        =\lambda_1\Oomega_1^a+\lambda_2\Oomega_2^a, & &
  \oomega_1^a=i(\mu {\overline\sigma}_1^a -{\overline\sigma}_2^a)
          =i(\lambda_1\Oomega_1^a-\lambda_2\Oomega_2^a),
  \nonumber\\
  \oomega_2^a=\mu {\overline\sigma}_2^a-{\overline\sigma}_1^a
          =\lambda_1\Oomega_2^a-\lambda_2\Oomega_1^a, & &
  \oomega_3^a=i({\overline\sigma}_1^a+\mu {\overline\sigma}_2^a)
          =i(\lambda_1\Oomega_2^a+\lambda_2\Oomega_1^a).
  \end{eqnarray}

  The differential $dp$ from the tangent
  sheaf to the pull-back of the
  tangent sheaf $p^*\co (2)$ on the projective line $\cp^1$
  is the  twisted 1-form
  $dp=\lambda_2d\lambda_1-\lambda_1d\lambda_2$ \cite[Example 13.83]{Besse}.
  Its kernel $\cd_Z$ is isomorphic to $\lie{t}^{1,0}\otimes \co (1)$.
  Applying Lemma \ref{sections} to $\ell=1$, we have
  \begin{lemma}\label{torus moduli}
  There is a natural isomorphism
  $
  H^k(Z, \cd_Z)=\lie{t}^{1,0}\otimes
   \lie{t}^{*(0,k)}\otimes S^{k+1}\bC^{2}.
  $
  \end{lemma}

  We seek local as well as global representations of elements in these cohomology
  spaces. On $U_2$ the product coordinates $\{z_1^a, z_2^a, \mu\}$ are not holomorphic.
  The holomorphic coordinates are
  \begin{equation}\label{holomorphic coordinate}
  w_1^a=\mu z_1^a-\zbar^a_2,
  \hspace{.2in}
  w_2^a=\mu z_2^a+\zbar^a_1,
  \hspace{.2in}
  \zeta=\mu.
  \end{equation}
  The inverse coordinate change is
  \begin{equation}
  z_1^a=\frac{1}{1+|\zeta |^2}\left( {\overline\zeta} w_1^a+\wbar^a_2\right),
  \hspace{.2in}
  z_2^a=\frac{1}{1+|\zeta |^2}\left( -\wbar^a_1+{\overline\zeta} w_2^a\right),
  \hspace{.2in}
  \mu=\zeta.
  \end{equation}
  When one changes coordinates from $\lambda_2\neq 0$ to
  $\lambda_1\neq 0$, $\mu{\tilde w}^a_j=w^a_j.$
  Therefore,
  \begin{equation}\label{vaj}
  V^a_j =\frac{1}{\lambda_2}\pdd{w^a_j}
  =\frac{1}{\lambda_1}\pdd{{\tilde w}^a_j}
  \end{equation}
  is a globally defined section of the tangent bundle.

  We are now able to describe a basis for the space $H^0(Z, \cd_Z)$.
  The dual of the (1,0)-form
  $\omega_k^a$ in (\ref{twisted forms}) is
  \begin{equation}\label{twisted fields}
  W_k^a=\frac{1}{2}\left(I_k\frac{\partial}{\partial x_{2a-1}}
      -iI_{\vec{a}}I_k\frac{\partial}{\partial x_{2a-1}}\right).
  \end{equation}
  These vector fields can also be identified as
  \begin{eqnarray}
  W_0^a=\lambda_1V_1^a+\lambda_2V_2^a,
  & &
  W_1^a=i(\lambda_1V_1^a-\lambda_2V_2^a),
  \nonumber\\
  W_2^a=\lambda_1V_2^a-\lambda_2V_1^a,
   & &
  W_3^a=i(\lambda_1V_2^a+\lambda_2V_1^a).
  \end{eqnarray}
  It follows that the elements
  $\lambda_1^{1-\ell}\lambda_2^\ell V_i^a\otimes\Omega_j^b$, with
  $0\leq \ell\leq 1$, $1\leq i,j\leq 2$, $1\leq a,b\leq m$
  form a basis for $H^1(Z, \cd_Z)$.
  Similarly,  the twisted vector-valued k-forms
  \begin{equation}
  \lambda_1^{k+1-\ell}\lambda_2^\ell V_i^a\otimes
  \Oomega_1^{a_1}\wedge\cdots\wedge\Oomega_1^{a_i}
        \wedge\Oomega_2^{b_1}\wedge\cdots\wedge\Oomega_2^{b_j}
  \end{equation}
   with $0\leq \ell\leq k+1$
  form a basis for $H^{k}(Z, \cd_Z)$.

  \subsection{The Twistor space of $X=\hei\times T^3$}
  From
  the definition of the complex structure on twistor
  spaces, the projection from $X$ to $T^{4m}$ induces a natural holomorphic
  projection $\Psi$ from the twistor space $W$ of $X$ onto the twistor
  space $Z$ of $T^{4m}$. Moreover, $p\circ \Psi=p$. We make use of these projections
  and the related spectral sequences to calculate cohomology on the twistor space $W$.

  \begin{lemma}\label{cohomology}
   For all $\ell\geq -1$,
  $H^k(W, p^*\co (\ell))
  = (\lie{h}_{4m+1}\oplus\lie{t}_3)^{*(0, k)}\otimes S^{k+\ell}\bC^{2}.$
  \end{lemma}
  \bproof
  The fibers of the projection $p$ from $W$ to $\cp^1$ are isomorphic to
  $\check{H}\times \U(1)\times T^2$ where $\check{H}\times U(1)$ is
  a complex $2m+1$-dimensional manifold \cite{GPP} and
  $T^2$ is a complex elliptic curve whose fundamental domain is
  a square. Since $H^1(T^2, \co)$ is one-dimensional, we
  combine \cite[Lemma 2]{GPP} with Kunneth formula to find that
  \begin{equation}
  H^k(p^{-1}(\vec{a}), \co)
  =\wedge^k H^1(p^{-1}(\vec{a}), \co).
  \end{equation}

  Let $\gamma$ be the dual of the central vector field $Z$ on the
  Heisenberg group. Denote
  the dual of $E_1, E_2$, and $ E_3$ by $\delta_1, \delta_2$ and $\delta_3$.
  As in (\ref{twisted forms}), define
  $\oomega_{k}^{m+1}=I_k\gamma-iI_{\vec{a}}I_k\gamma.$
  Given the stereographic projection (\ref{stereographic}), we deduce that
  when $\lambda_2\neq 0$,
  \begin{eqnarray}
  \oomega_0^{m+1}=\mu {\overline\sigma}_1^{m+1}+{\overline\sigma}_2^{m+1},
  & &
  \oomega_1^{m+1}=i(\mu {\overline\sigma}_1^{m+1}-{\overline\sigma}_2^{m+1}),
  \nonumber\\
  \oomega_2^{m+1}=\mu {\overline\sigma}_2^{m+1}-{\overline\sigma}_1^{m+1},
  & &
  \oomega_3^{m+1}=i({\overline\sigma}_1^{m+1}+\mu {\overline\sigma}_2^{m+1}),
  \end{eqnarray}
  where
  \begin{equation}\label{sigma(m+1)}
  \osigma_1^{m+1}=
  \frac{\mubar (\gamma-i\delta_1)-(\delta_2+i\delta_3)}{1+|\mu |^2},
  \hspace{.2in}
  \osigma_2^{m+1}
  =\frac{(\gamma+i\delta_1)+{\overline\mu} (\delta_2-i\delta_3)}{1+| \mu |^2}.
  \end{equation}
  It is important to note that $d\gamma$ is type (1,1)
  with respect to any complex structure in the given hypercomplex structure.
  In fact, the structural equation yields
  \begin{equation}
  d\gamma=4\sum_{a=1}^m(dx_{2a-1}\wedge dy_{2a-1}+dx_{2a}\wedge dy_{2a}).
  \end{equation}
  It follows that
  \begin{eqnarray}
  d{\overline\sigma}_1^{m+1}
  &=&2{\overline\mu}\sum_{a=1}^m(\sigma_1^a\wedge{\overline\sigma}_2^a
       +{\overline\sigma}_1^a\wedge\sigma_2^a)
       +\frac{d{\overline\mu}\wedge\sigma_2^{m+1}
          -{\overline\mu}d\mu\wedge{\overline\sigma}_1^{m+1}}
            {1+|\mu|^2};\label{dsigma 1}
  \\
  d{\overline\sigma}_2^{m+1}
  &=&2\sum_{a=1}^m(\sigma_1^a\wedge{\overline\sigma}_2^a
       +{\overline\sigma}_1^a\wedge\sigma_2^a)
       +\frac{-{\overline\mu}d\mu\wedge{\overline\sigma}_2^{m+1}
          -d{\overline\mu}\wedge{\sigma}_1^{m+1}}
            {1+|\mu|^2}. \label{dsigma 2}
  \end{eqnarray}
  In particular,
  $\osigma_1^{m+1}$ and $\osigma_2^{m+1}$ are holomorphic.
  It shows that the (0,1)-forms $\oomega_k^{m+1}$ are holomorphic
  on the twistor space $W$.
  These holomorphic (0,1)-forms  are the analogies of
  $\osigma_1^a$ and $\osigma_2^a$
  on the twistor space $Z$ defined in (\ref{fourteen}). These 1-forms on
  $Z$ are pulled back by $\Psi$ to holomorphic (0,1)-forms on $W$.
  Therefore,
  $\osigma_i^\alpha$,  with $i=1,2$ and $1\leq \alpha\leq m+1$
  form a basis for
  $H^1(p^{-1}(\vec{a}), \co)$ when $\lambda_2\neq 0$. Due to
  the homogeneity of these forms,
  \begin{equation}
  R^q p_*\co_W=(\lie{h}_{4m+1}\oplus \lie{t}_3)^{*(0,q)}\otimes \co(q).
  \end{equation}

  Consider the Leray spectral
  sequence $E_*^{p,q}(p^*\co (\ell))$. Due to the projection formula,
  \begin{eqnarray}
  & &E_2^{p,q}(p^*\co (\ell)) = H^p(\cp^1, R^q p_*p^*\co (\ell))
    =H^p(\cp^1, \co (\ell)\otimes R^q p_*\co_W)
  \nonumber\\
          &=& (\lie{h}_{4m+1}\oplus \lie{t}_3)^{*(0,q)}\otimes
         H^p(\cp^1, \co (q+\ell )).
  \end{eqnarray}
  As every element in this space is represented by global holomorphic
  forms on the twistor space $W$ with
  coefficients in $p^*\co (\ell) $,
  $d_2=0$. It follows that for $\ell\geq -1$,
  \begin{eqnarray*}
  H^k(W, p^*\co (\ell) ) &=& \oplus_{p+q=k}E_2^{p,q}(p^*\co (\ell))
        =E_2^{0,k}(p^*\co (\ell))
  \\
  &=& (\lie{h}_{4m+1}\oplus\lie{t}_3)^{*(0, k)}\otimes S^{k+\ell}\bC^{2}.
  \end{eqnarray*}
  \eproof

  Since the kernel $\cv$ of the differential $d\Psi$ restricted to
  $\cd_W$ satisfies $\cv=(\lie{c}\oplus \lie{t}_3)^{1,0}\otimes \co (1)$, and
  since $\Psi^*\cd_Z=\lie{t}_{4m}^{1,0}\otimes \co (1)$, the next lemma follows.

  \begin{lemma}\label{number5} Let $\cv$ be the kernel of the differential $d\Psi$ restricted to
  $\cd_W$. Then
  \begin{eqnarray*}
  H^k(W, \cv) &=& (\lie{c}\oplus \lie{t}_3)^{1,0}\otimes
     (\lie{h}_{4m+1}\oplus\lie{t}_3)^{*(0, k)}\otimes S^{k+1}\bC^{2},
  \\
  H^k(W, \Psi^*\cd_Z) &=&
         \lie{t}_{4m}^{1,0}\otimes
      (\lie{h}_{4m+1}\oplus\lie{t}_3)^{*(0, k)}\otimes S^{k+1}\bC^{2}.
  \end{eqnarray*}
  \end{lemma}

  This lemma implies that $H^0(W, \cv)=(\lie{c}\oplus \lie{t}_3)^{1,0}\otimes \bC^{2}$.
  This space is linearly spanned by the infinitesimal
  hypercomplex transformations generated
  by the center of the algebra $\lie{h}_{4m+1}\oplus \lie{t}_3$. It
  is a linear span of  $W_k:=\frac12 (I_kZ-iI_{\vec{a}}I_kZ)$.

  On $W$, a local holomorphic frame for the bundle of (0,1)-forms consists of
  $\{\osigma_i^\alpha, d{\overline{\mu}}\}$.
  The dual smooth (0,1)-vectors are $\partial/\partial\mubar$ and
  \begin{eqnarray}
  \dbar_1^a &=& \frac12\{\mu (X_{2a-1}+iX_{2a})-(Y_{2a-1}-iY_{2a}) \},
  \nonumber\\
  \dbar_2^a &=& \frac12\{ (X_{2a-1}-iX_{2a})+\mu (Y_{2a-1}+iY_{2a}) \},
  \label{32}\\
  \dbar_1^{m+1} &=& \frac12\{\mu (Z+iE_1)-(E_2-iE_3) \},
  \dbar_2^{m+1} = \frac12\{  (Z-iE_1)+\mu (E_2+iE_3) \}.
  \label{34}
  \end{eqnarray}

  Now we examine the induced long exact sequence of
  \begin{equation}\label{exact*}
  0\to {\cal V} \to {\cal D}_W \to \Psi^*{\cal D}_Z\to 0
  \end{equation}
  to calculate the cohomology of ${\cal D}_W$.
  Note that elements in $H^0(W, \Psi^*{\cd}_Z)$ are
  linear combinations of the vector fields (\ref{twisted fields}).
  They have natural lifting to smooth sections of the
  tangent bundle on $W$, namely
  ${\tilde{W}}_a^k=\frac12 (I_kX_{2a-1}-iI_{\vec{a}}I_kX_{2a-1})$.
  Given the smooth local frame above,
  \begin{eqnarray}\label{the ws}
  {\tilde{W}}_0^a=\frac{1}{1+|\mu|^2}(\mu\partial_1^a+\partial_2^a),
  & &
  {\tilde{W}}_1^a=\frac{i}{1+|\mu|^2}(\mu\partial_1^a-\partial_2^a),
  \nonumber\\
  {\tilde{W}}_2^a=\frac{1}{1+|\mu|^2}(\mu\partial_2^a-\partial_1^a),
   & &
  {\tilde{W}}_3^a=\frac{i}{1+|\mu|^2}(\mu\partial_2^a+\partial_1^a).
  \end{eqnarray}
  Given the algebra structure of the Heisenberg algebra,
  \begin{equation}\label{commutator of dbar}
  [\dbar_i^a, \partial_i^b]=0,
   [\dbar_2^a, \partial_1^b]
  =-[\dbar_1^a, \partial_2^b]
  ={2\delta_{ab}}(1+|\mu|^2)Z.
  \hspace{.3in}
  [\partial^a_i, \partial^b_j]=0.
  \end{equation}
  In particular,
  \begin{equation}\label{gauduchon}
   [\dbar_2^a, \partial_1^a]^{1,0}
  =-[\dbar_1^a, \partial_2^a]^{1,0}
  =2(1+|\mu|^2)W_0.
  \end{equation}
  A computation using  (\ref{32}) and (\ref{the ws}) shows that
  \begin{eqnarray}
  [\pdd{\mubar}, {\tilde{W}}_0^a]
      = \frac{1}{(1+|\mu|^2)^2}\left( -\dbar_1^a+\mu\dbar_2^a\right),
  &\ &
  {[{\pdd{\mubar}}, {\tilde{W}}_1^a]}
     = \frac{i}{(1+|\mu|^2)^2}\left( \dbar_1^a+\mu\dbar_2^a\right),
  \\
  {[{\pdd{\mubar}}, {\tilde{W}}_2^a]}
     = \frac{-1}{(1+|\mu|^2)^2}\left( \mu\dbar_1^a+\dbar_2^a\right),
  &\ &
  {[\pdd{\mubar}, {\tilde{W}}_3^a]}
     = \frac{i}{(1+|\mu|^2)^2}\left( -\mu\dbar_1^a+\dbar_2^a\right).
  \end{eqnarray}
  In particular, $[\pdd{\mubar}, {\tilde{W}}_k^a]^{1,0}=0$ for all $a$ and
  $k$. Note also that $[\dbar_j^{m+1}, {\tilde{W}}_k^a]=0$ because
  $Z, E_1, E_2, E_3$ are in the center of $\lie{h}_{4m+1}\oplus \lie{t}_3$.
  It follows that the 0-th coboundary map for the induced exact sequence of
  (\ref{exact*}) can be calculated as follows.
  We use the Chern connection $\nabla$
  on the holomorphic tangent bundle on $W$ to define the $\dbar$-operator
  $\dbar^\nabla$. Then $\delta_0 (W_k^a)$ is represented uniquely by the cohomology
  class of $\dbar^\nabla {\tilde{W}}_k^a$.
  Due to an observation of Gauduchon \cite{Gauduchon},
  for any (0,1)-vector $X$, and (1,0)-vector $Y$,
  $\dbar^\nabla_X Y=[X, Y]^{1,0}$. The above discussion implies that
  \[
  \delta_0(W_k^a)
  = [\dbar_j^a, {\tilde{W}}_k^a]^{1,0}\otimes \osigma_j^a
  =[\dbar_1^a, {\tilde{W}}_k^a]^{1,0}\otimes \osigma_1^a
  +[\dbar_2^a, {\tilde{W}}_k^a]^{1,0}\otimes \osigma_2^a.
  \]
  With (\ref{the ws}), (\ref{commutator of dbar}) and (\ref{gauduchon}), we deduce that
  \begin{eqnarray*}
  \delta_0(W_0^a)&=&2W_0\otimes
    (\lambda_1\Oomega_2^a-\lambda_2 \Oomega_1^a),
  \delta_0(W_1^a)=2iW_0\otimes
    (\lambda_1\Oomega_2^a+\lambda_2 \Oomega_1^a),
  \\
  \delta_0(W_2^a)&=&-2W_0\otimes
    (\lambda_1\Oomega_1^a+\lambda_2 \Oomega_2^a),
  \delta_0(W_3^a)=-2iW_0\otimes
    (\lambda_1\Oomega_1^a-\lambda_2 \Oomega_2^a).
  \end{eqnarray*}
  Therefore, the coboundary map $\delta_0$
  is injective and
  $
  H^0(W, \cd_W) = H^0(W, \cv)$.

  The functions $w_1^{m+1}=\mu z_1-\zbar_2$ and
  $w_2^{m+1}=\mu z_2+\zbar_1$ are holomorphic on the open subset $U_2$ of
   the twistor space $W$. On $U_1$, define
  ${\tilde w}_j^{m+1}=\frac{1}{\mu}w_j^{m+1}$. As in (\ref{vaj}) we define
  $V_j:=V_j^{m+1}: =\frac{1}{\lambda_2}\pdd{w_j^{m+1}}
  =\frac{1}{\lambda_1}\pdd{{\tilde w}_j^{m+1}}$.
  In this basis,
  the image of $H^0(W, \Psi^*\cd_Z)$ in $H^1(W, \cv)$ is spanned by
  $\lambda_i (\lambda_1V_1+\lambda_2V_2)\otimes \Oomega_j^a$.
  By Lemma \ref{number5},
  \[
  H^1(W, {\cal V}) =
   (\lie{c}\oplus\lie{t}_3)^{1,0}
        \otimes (\lie{c}\oplus \lie{t}_3)^{*(0,1)}\otimes S^2\bC^2
     \oplus (\lie{c}\oplus\lie{t}_3)^{1,0}
        \otimes \lie{t}_{4m}^{*(0,1)}\otimes S^2\bC^2.
  \]
  Since the image of $\delta_0$ is contained in the second summand, the cokernel
  of $\delta_0$ is the direct sum of the following spaces.
  \begin{eqnarray}
  coker{\delta_0}'': &=&(\lie{c}\oplus\lie{t}_3)^{1,0}
        \otimes (\lie{c}\oplus \lie{t}_3)^{*(0,1)}\otimes S^2\bC^2
  = \mbox{span}\{\lambda_1^{k}\lambda_2^{2-k}V_i\Oomega_j^{m+1}\},
  \nonumber \\
  coker{\delta_0}': &=&
     \mbox{span}\{ (\lambda_1V_1-\lambda_2V_2)\lambda_i\Oomega_j^a,
               \lambda_2^2V_1\Oomega_j^a, \lambda_1^2V_2\Oomega_j^a\},
  \end{eqnarray}
  where $\Oomega_j^{m+1}=\frac{1}{\lambda_2}\osigma_j^{m+1}$.

  A basis for
  $H^1(W, \Psi^*\cd_Z)$ consists of
  $\lambda_1^{2-k}\lambda_2^k V_i^a\Oomega_j^\beta$,  where
  $0\leq k\leq 2$, $1\leq i, j\leq 2$, $1\leq a\leq m$, $1\leq \beta\leq m+1$.
  Since $\Oomega_j^\beta$ is holomorphic, the coboundary map $\delta_1$ from
  the first cohomology to the second has the following property.
  \[
  \delta_1(V_1^a\Oomega_j^\beta)
  =2(\lambda_1V_1+\lambda_2V_2)\otimes \Oomega_2^a\wedge
     \Oomega_j^{\beta},
  \
  \delta_1(V_2^a\Oomega_j^\beta)
  =-2(\lambda_1V_1+\lambda_2V_2)\otimes\Oomega_1^a\wedge
      \Oomega_j^\beta.
  \]
  Therefore, the kernel of $\delta_1$ is
  \begin{equation}\label{kernel of d1}
  \ker{\delta_1}=\mbox{span}\{\lambda_1^k\lambda_2^{2-k}\}\otimes
  \mbox{span}\{V_1^a\Oomega_2^b+V_1^b\Oomega_2^a,
            V_2^a\Oomega_1^b+V_2^b\Oomega_1^a,
    V_1^a\Oomega_1^b-V_2^b\Oomega_2^a\}.
  \end{equation}
  The induced exact sequence yields
  \begin{eqnarray}
   H^1(W, \cd_W) &=&coker{\delta_0}''\oplus coker{\delta_0}'\oplus ker{\delta_1}
  \nonumber \\
  &=&
  \{\lambda_1^{k}\lambda_2^{2-k}V_i\Oomega_j^{m+1}\}\oplus
  \{ (\lambda_1V_1-\lambda_2V_2)\lambda_i\Oomega_j^a,
               \lambda_2^2V_1\Oomega_j^a, \lambda_1^2V_2\Oomega_j^a\}
  \nonumber\\
  &\oplus &\{\lambda_1^k\lambda_2^{2-k}\}\otimes
  \{V_1^a\Oomega_2^b+V_1^b\Oomega_2^a,
            V_2^a\Oomega_1^b+V_2^b\Oomega_1^a,
    V_1^a\Oomega_1^b-V_2^b\Oomega_2^a\}.
  \end{eqnarray}

  Here we see that $coker {\delta_0}''$ is a $12$-dimensional space. It is generated
  by vertical tangent vectors and 1-forms with respect to the projection $\Psi$.
  The dimension of $coker{\delta_0}'$ is equal to $8m$. The contribution
  from deformation of the basis of the projection $\Psi$ is in $ker{\delta_1}$.
  The dimension of this space is equal to $3m(2m+1)$. This concludes the
  proof of Theorem \ref{theorem1}.

  \section{Deformation Theory}

  The obstruction space to hypercomplex deformation is $H^2(W, \cd_W)$. Results
  in the last section demonstrate that it does not vanish, we approach
  integrability of parameters in $H^1(W, \cd_W)$ by constructing
  a one-parameter deformation for any given tangent. It is done
  through Kodaira-Spencer's
  method of constructing a convergent sequence with coefficients in
  vector-valued 1-forms.
  During this construction,
  we do not control the vector-valued 1-forms
  in such a way that the vector part is tangent to the distribution $\cd_W$.
  Therefore, the deformation may a priori depart from hypercomplex structures. We
  tackle this problem at the end of this section.

  \subsection{Integrability}\label{integrable}
  We have seen that the virtual parameter space $H^1(W, \cd)$ is a vector subspace of
  the following space.
  \begin{eqnarray}\label{spaceE}
  \ce &= &H^1(W, \cv)\oplus ker{\delta_1}
  \nonumber
  \\
   &=& (\lie{c}\oplus \lie{t}_3)^{1,0}
     \otimes(\lie{h}_{4m+1}\oplus\lie{t}_3)^{*(0,1)}\otimes S^2\bC^2
  \oplus ker{\delta_1}
  \nonumber \\
  &=& \langle \lambda^k_1\lambda_2^{2-k}\rangle \otimes
     \left( \langle V_i^{m+1} \Oomega_j^\beta \rangle
       \oplus
  \langle V_1^a\Oomega_2^b+V_1^b\Oomega_2^a,
            V_2^a\Oomega_1^b+V_2^b\Oomega_1^a,
    V_1^a\Oomega_1^b-V_2^b\Oomega_2^a
  \rangle
   \right).
  \end{eqnarray}
  We denote these elements
  by $\Up_\tau$, where $1\leq \tau
  \leq D=\dim \ce$.
  If  $\Omega_1=\omega\otimes V$ and $\Omega_2=\omega'
  \otimes V'$ are vector-valued 1-forms,
   the Dolbeault representative for the Nijenhuis bracket
  $\{\omega\otimes V, \omega'\otimes V'\}$ is \cite{Nijenhuis, N2}
  \begin{equation}\label{Kuranishi bracket}
  \omega'\wedge L_{V'}\omega\otimes V
   +\omega\wedge L_{V}\omega'\otimes V'
   +\omega\wedge\omega' \otimes [V, V'].
  \end{equation}
  The aim of this section is to prove the following proposition.
  \begin{proposition}\label{expansion}
  Let $\Gamma^0$ be the complex vector space of complex-valued
  smooth functions on $S^2$. For any $\phi_1$ and $\phi_2$ in
  $\Gamma^0\otimes {\ce}$,
  there exist smooth functions $g_\tau$ on the 2-sphere such that
  $\{\phi_1, \phi_2\}=\sum_{\tau}\dbar\left( g_\tau\Up_\tau\right).$
  \end{proposition}
  The following
  functions on 2-sphere  $S^2$ is extended to the twistor space $W$.
  \begin{equation}
  f_1= \frac{\mu}{1+|\mu|^2},
  \hspace{.2in}
  f_2= \frac{{\overline\mu}}{1+|\mu|^2},
  \hspace{.2in}
  f_3= \frac{1}{1+|\mu|^2}.
  \end{equation}
  We have
  \begin{equation}
  \dbar f_1=-\frac{\mu^2 d{\overline\mu}}{(1+|\mu|^2)^2},
  \hspace{.2in}
  \dbar f_2 = \frac{ d{\overline\mu}}{(1+|\mu|^2)^2},
  \hspace{.2in}
  \dbar f_3=-\frac{\mu d{\overline\mu}}{(1+|\mu|^2)^2}.
  \end{equation}
  In subsequent computation, we make use of these three functions to prove
  the Proposition \ref{expansion}.

  Let  $\epsilon_{11}=\epsilon_{22}=0$ and $\epsilon_{12}=-\epsilon_{21}=1$.
  It follows  from (\ref{32}) and (\ref{34}) that
  \begin{equation}
  [\partial_i^\alpha, \frac{\partial}{\partial{\overline\mu}}]^{0,1}
  =-\frac{\epsilon_{ij}}{1+|\mu|^2}\partial_j^\alpha,
  \mbox{ and }
  L_{\partial_i^\alpha}\osigma_j^{\beta}=\delta_{\alpha\beta}\epsilon_{ij}
    \frac{d{\overline\mu}}{1+|\mu|^2}
  \end{equation}
  except when $\alpha=a, \beta=m+1$.
  Therefore, in the case when $\alpha=\alpha'=m+1$ or when
  $1\leq \beta, \beta'\leq m$, with (\ref{commutator of dbar}) we deduce that
  \begin{eqnarray*}
  & &\{ \lambda_1^k\lambda_2^{2-k}V_i^\alpha\Oomega_j^\beta,
      \lambda_1^{k'}\lambda_2^{2-k'}V_{i'}^{\alpha'}\Oomega_{j'}^{\beta'}\}
  \\
  &=& \{ \frac{\mu^k}{1+|\mu|^2}\partial_i^\alpha\osigma_j^\beta,
   \frac{\mu^{k'}}{1+|\mu|^2}\partial_{i'}^{\alpha'}\osigma_{j'}^{\beta'}\}
  = \frac{\mu^{k+k'}}{(1+|\mu|^2)^2}
    \{\partial_i^\alpha\osigma_j^\beta,
   \partial_{i'}^{\alpha'}\osigma_{j'}^{\beta'}\}
  \\
  &=& \frac{\mu^{k+k'}}{(1+|\mu|^2)^2}
   \left(
  \partial_i^\alpha
  \osigma_{j'}^{\beta'}\wedge L_{\partial_{i'}^{\alpha'}}
   (\osigma_j^\beta)
  +
  \partial_{i'}^{\alpha'}
  \osigma_j^\beta\wedge L_{\partial_{i}^{\alpha}}
   (\osigma_{j'}^{\beta'})
  +[\partial_i^\alpha, \partial_{i'}^{\alpha'}]
  \osigma_j^\beta\wedge \osigma_{j'}^{\beta'}
   \right)
  \\
  &=& \frac{\mu^{k+k'}}{(1+|\mu|^2)^2}
   \left(
  \partial_i^\alpha
  \osigma_{j'}^{\beta'}\wedge L_{\partial_{i'}^{\alpha'}}
   (\osigma_j^\beta)
  +
  \partial_{i'}^{\alpha'}
  \osigma_j^\beta\wedge L_{\partial_{i}^{\alpha}}
   (\osigma_{j'}^{\beta'})
   \right)
  \\
  &=& -\frac{\mu^{k+k'}}{(1+|\mu|^2)^3}
   d{\overline\mu}\wedge \left(
  \delta_{\alpha'\beta} \epsilon_{i'j}
  \partial_i^\alpha \osigma_{j'}^{\beta'}
  +
  \delta_{\alpha\beta'} \epsilon_{ij'}
  \partial_{i'}^{\alpha'}
  \osigma_j^\beta
  \right)
  \\
  &=&
  -\frac{\mu^{k'}d{\overline\mu}}{(1+|\mu|^2)^2}
   \wedge \left(
  \frac{\delta_{\alpha'\beta} \epsilon_{i'j}
  \mu^{k}}{1+|\mu|^2}
  \partial_i^\alpha \osigma_{j'}^{\beta'}
  \right)
  -\frac{\mu^{k'}d{\overline\mu}}{(1+|\mu|^2)^2}
   \wedge \left(
  \frac{\delta_{\alpha\beta'} \epsilon_{ij'}
  \mu^{k}}{1+|\mu|^2}
  \partial_{i'}^{\alpha'}
  \osigma_j^\beta
  \right)
  \\
  &=&
  -\frac{\mu^{k'}d{\overline\mu}}{(1+|\mu|^2)^2}
   \wedge \left(
  \delta_{\alpha'\beta} \epsilon_{i'j}
  \lambda_1^k\lambda_2^{2-k}V_i^\alpha\Oomega_{j'}^{\beta'}\right)
  -\frac{\mu^{k'}d{\overline\mu}}{(1+|\mu|^2)^2}
   \wedge \left(
  \delta_{\alpha\beta'} \epsilon_{ij'}
  \lambda_1^{k}\lambda_2^{2-k}V_{i'}^{\alpha'}\Oomega_j^\beta
  \right)
  \\
  &=& \dbar
  \left( -f \delta_{\alpha'\beta} \epsilon_{i'j}
  \lambda_1^k\lambda_2^{2-k}V_i^\alpha\Oomega_{j'}^{\beta'}
  -f \delta_{\alpha\beta'} \epsilon_{ij'}
  \lambda_1^{k}\lambda_2^{2-k}V_{i'}^{\alpha'}\Oomega_j^\beta
  \right)
  \\
  &=& \dbar
  \left( (-f\lambda_1^k\lambda_2^{2-k})( \delta_{\alpha'\beta} \epsilon_{i'j}
  V_i^\alpha\Oomega_{j'}^{\beta'}
  +\delta_{\alpha\beta'} \epsilon_{ij'}
  V_{i'}^{\alpha'}\Oomega_j^\beta)
  \right).
  \end{eqnarray*}
  where $f$ is determined by $k'$.

  \begin{lemma}For any $\phi_1$ and $\phi_2$ in
  $H^1(W, \cv)$,
  there exist smooth functions $g_\tau$ on the 2-sphere and
  $\phi_\tau\in H^1(W, \cv)$ such that
  $\{\phi_1, \phi_2\}=\sum_{\tau}\dbar\left( g_\tau\phi_\tau\right).$
  \end{lemma}
  \bproof By Lemma \ref{number5},  elements in $H^1(W, \cv)$ have the form
  $\lambda_1^k\lambda_2^{2-k}V_i^{m+1}\Oomega_j^\beta$ for some $k, j$ and
  $\beta$.
   From the above computation,
  \begin{eqnarray*}
  & &\{ \lambda_1^k\lambda_2^{2-k}V_i^{m+1}\Oomega_j^\beta,
      \lambda_1^{k'}\lambda_2^{2-k'}V_{i'}^{m+1}\Oomega_{j'}^{\beta'}\}
  \\
  &=& \dbar
  \left( (-f\lambda_1^k\lambda_2^{2-k})( \delta_{m+1,\beta} \epsilon_{i'j}
  V_i^{m+1}\Oomega_{j'}^{\beta'}
  +\delta_{m+1,\beta'} \epsilon_{ij'}
  V_{i'}^{m+1}\Oomega_j^\beta)
  \right).
  \end{eqnarray*}
  This is an element in the image of $\Gamma^0\otimes H^1(W, \cv)$ via $\dbar$.
  \eproof

  \begin{lemma}For any $\phi_1$ and $\phi_2$ in
  $\ker \delta_1$,
  there exist smooth functions $g_\tau$ on the 2-sphere and
  $\phi_\tau\in \ker\delta_1$ such that
  $\{\phi_1, \phi_2\}=\sum_{\tau}\dbar\left( g_\tau\phi_\tau\right).$
  \end{lemma}
  \bproof
  We do it through a case-by-case computation. In view of formula
  (\ref{kernel of d1}), there are six cases.
  \begin{eqnarray*}
  & &\{\lambda_1^k\lambda_2^{2-k}
      (V_1^a\Oomega_2^b+V_1^b\Oomega_2^a),
         \lambda_1^{k'}\lambda_2^{2-k'}
     (V_1^{a'}\Oomega_2^{b'}+V_1^{b'}\Oomega_2^{a'})\}
  \\
  &=&
  \dbar
  \left( -f\lambda_1^k\lambda_2^{2-k}
  ( \delta_{a'b}  V_1^a\Oomega_2^{b'}
  +\delta_{ab'}   V_1^{a'}\Oomega_2^b)
  \right)
  +
  \dbar
  \left( -f\lambda_1^k\lambda_2^{2-k}
  ( \delta_{b'b}  V_1^a\Oomega_2^{a'}
  +\delta_{aa'}   V_1^{b'}\Oomega_2^b)
  \right)
  \\
  & &
  +
  \dbar
  \left( -f\lambda_1^k\lambda_2^{2-k}
  ( \delta_{a'a}  V_1^b\Oomega_2^{b'}
  +\delta_{bb'}   V_1^{a'}\Oomega_2^a)
  \right)
  +
  \dbar
  \left( -f\lambda_1^k\lambda_2^{2-k}
  ( \delta_{b'a}  V_1^b\Oomega_2^{a'}
  +\delta_{ba'}   V_1^{b'}\Oomega_2^a)
  \right)
  \\
  &=&
  \dbar
  \left( -f\lambda_1^k\lambda_2^{2-k} \delta_{a'b}
  ( V_1^a\Oomega_2^{b'}+ V_1^{b'}\Oomega_2^a)\right)
  +
  \dbar
  \left( -f\lambda_1^k\lambda_2^{2-k} \delta_{ab'}
  ( V_1^{a'}\Oomega_2^{b}+ V_1^{b}\Oomega_2^{a'})\right)
  \\
  & &+\dbar
  \left( -f\lambda_1^k\lambda_2^{2-k} \delta_{bb'}
  ( V_1^a\Oomega_2^{a'}+ V_1^{a'}\Oomega_2^a)\right)
  +
  \dbar
  \left( -f\lambda_1^k\lambda_2^{2-k} \delta_{aa'}
  ( V_1^{b}\Oomega_2^{b'}+ V_1^{b'}\Oomega_2^{b})\right).
  \end{eqnarray*}
  This element is in the image of $\Gamma^0\otimes \ker\delta_1$ via $\dbar$.
  Similarly, when we interchange the indices 1 and 2, we have
  \begin{eqnarray*}
  & &\{\lambda_1^k\lambda_2^{2-k}
      (V_2^a\Oomega_1^b+V_2^b\Oomega_1^a),
         \lambda_1^{k'}\lambda_2^{2-k'}
     (V_2^{a'}\Oomega_1^{b'}+V_2^{b'}\Oomega_1^{a'})\}
  \\
  &=&
  \dbar
  \left( f\lambda_1^k\lambda_2^{2-k} \delta_{a'b}
  ( V_2^a\Oomega_1^{b'}+ V_2^{b'}\Oomega_1^a)\right)
  +
  \dbar
  \left( f\lambda_1^k\lambda_2^{2-k} \delta_{ab'}
  ( V_2^{a'}\Oomega_1^{b}+ V_2^{b}\Oomega_1^{a'})\right)
  \\
  & &+\dbar
  \left( f\lambda_1^k\lambda_2^{2-k} \delta_{bb'}
  ( V_2^a\Oomega_1^{a'}+ V_2^{a'}\Oomega_1^a)\right)
  +
  \dbar
  \left( f\lambda_1^k\lambda_2^{2-k} \delta_{aa'}
  ( V_2^{b}\Oomega_1^{b'}+ V_2^{b'}\Oomega_1^{b})\right).
  \end{eqnarray*}
  This element is in the image of $\Gamma^0\otimes \ker\delta_1$ via $\dbar$. Next,
  \begin{eqnarray*}
  & &\{\lambda_1^k\lambda_2^{2-k}
      (V_1^a\Oomega_1^b-V_2^b\Oomega_2^a),
         \lambda_1^{k'}\lambda_2^{2-k'}
     (V_1^{a'}\Oomega_1^{b'}-V_2^{b'}\Oomega_2^{a'})\}
  \\
  &=&-\dbar \left(
     f \lambda_1^k\lambda_2^{2-k}
       (-\delta_{b'b}V_1^a\Oomega_2^{a'}+\delta_{aa'}V_2^{b'}\Oomega_1^b)
      \right)
  -\dbar \left(
     f \lambda_1^k\lambda_2^{2-k}
       (\delta_{a'a}V_2^b\Oomega_1^{b'}-\delta_{bb'}V_1^{a'}\Oomega_2^a)
      \right)
  \\
  &=&
  \dbar \left(
     f \lambda_1^k\lambda_2^{2-k}
     \delta_{b'b}(V_1^a\Oomega_2^{a'}+V_1^{a'}\Oomega_2^a)
      \right)
  -\dbar \left(
     f \lambda_1^k\lambda_2^{2-k}
      \delta_{a'a}(V_2^b\Oomega_1^{b'}+V_2^{b'}\Oomega_1^b)
      \right)
  \end{eqnarray*}
  This element is in the image of $\Gamma^0\otimes \ker\delta_1$ via $\dbar$.
  On the other hand,
  \[
  \{\lambda_1^k\lambda_2^{2-k}
      (V_1^a\Oomega_2^b+V_1^b\Oomega_2^a),
         \lambda_1^{k'}\lambda_2^{2-k'}
     (V_2^{a'}\Oomega_1^{b'}+V_2^{b'}\Oomega_1^{a'})\}=0.
  \]
  \begin{eqnarray*}
  & &
  \{\lambda_1^k\lambda_2^{2-k}
      (V_1^a\Oomega_2^b+V_1^b\Oomega_2^a),
         \lambda_1^{k'}\lambda_2^{2-k'}
     (V_1^{a'}\Oomega_1^{b'}-V_2^{b'}\Oomega_2^{a'})\}
  \\
  &=&
  \dbar(f\lambda_1^k\lambda_2^{2-k}\delta_{a'b}V_1^a\Oomega_1^{b'})
  -\dbar(f\lambda_1^k\lambda_2^{2-k}\delta_{aa'}V_2^{b'}\Oomega_2^{b})
  \\
  & &
  +\dbar(f\lambda_1^k\lambda_2^{2-k}\delta_{a'a}V_1^b\Oomega_1^{b'})
  -\dbar(f\lambda_1^k\lambda_2^{2-k}\delta_{ba'}V_2^{b'}\Oomega_2^{a})
  \\
  &=&
  \dbar(f\lambda_1^k\lambda_2^{2-k}\delta_{a'b}
         (V_1^a\Oomega_1^{b'}-V_2^{b'}\Oomega_2^{a}))
  +
  \dbar(f\lambda_1^k\lambda_2^{2-k}\delta_{a'a}
      (V_1^b\Oomega_1^{b'}-V_2^{b'}\Oomega_2^{b})).
  \end{eqnarray*}
  This element is in the image of $\Gamma^0\otimes \ker\delta_1$ via $\dbar$.
  Finally,
  \begin{eqnarray*}
  & &
  \{\lambda_1^k\lambda_2^{2-k}
      (V_2^a\Oomega_1^b+V_2^b\Oomega_1^a),
         \lambda_1^{k'}\lambda_2^{2-k'}
     (V_1^{a'}\Oomega_1^{b'}-V_2^{b'}\Oomega_2^{a'})\}
  \\
  &=&
  -\dbar(f\lambda_1^k\lambda_2^{2-k}\delta_{ab'}V_1^{a'}\Oomega_1^{b})
  +\dbar(f\lambda_1^k\lambda_2^{2-k}\delta_{b'b}V_2^{a}\Oomega_2^{a'})
  \\
  & &
  -\dbar(f\lambda_1^k\lambda_2^{2-k}\delta_{bb'}V_1^{a'}\Oomega_1^a)
  +\dbar(f\lambda_1^k\lambda_2^{2-k}\delta_{b'a}V_2^{b}\Oomega_2^{a'})
  \\
  &=&
  -\dbar(f\lambda_1^k\lambda_2^{2-k}\delta_{ab'}
         (V_1^{a'}\Oomega_1^{b}-V_2^{b}\Oomega_2^{a'}))
  -
  \dbar(f\lambda_1^k\lambda_2^{2-k}\delta_{bb'}
      (V_1^{a'}\Oomega_1^{a}-V_2^{a}\Oomega_2^{a'})).
  \end{eqnarray*}
  This element is in the image of $\Gamma^0\otimes \ker\delta_1$ via $\dbar$.
  \eproof

  To consider the Nijenhuis bracket between elements in $H^1(W, \cv)$ and
  elements in $\ker\delta_1$, we recall that $H^1(W, \cv)$ has two types of elements:
   $\lambda_1^k\lambda_2^{2-k}V_i^{m+1}\Oomega_j^{m+1}$ as elements in $coker\delta_0"$
  and $\lambda_1^k\lambda_2^{2-k}V_i^{m+1}\Oomega_j^b$ as elements in the
   direct summand complement in $H^1(W, \cv)$. We separate the computation into
  the next two lemmas.
  This computation  involves the algebraic structure of
  Heisenberg group because by (\ref{dsigma 1}) and (\ref{dsigma 2}),
  \begin{equation}
  L_{\partial_1^a}{\osigma}_1^{m+1}=2\mubar\osigma_2^a,
  \quad
  L_{\partial_2^a}{\osigma}_1^{m+1}=-2\mubar\osigma_1^a,
  \quad
  L_{\partial_1^a}{\osigma}_2^{m+1}=2\osigma_2^a,
  \quad
  L_{\partial_2^a}{\osigma}_2^{m+1}=-2\osigma_1^a.
  \end{equation}

  \begin{lemma}For any $\phi_1$ in $ker\delta_1$ and
  $\phi_2=\lambda_1^{k'}\lambda_2^{2-k'}V_{i'}^{m+1}\Oomega_{j'}^{b'}$,
  there exist smooth functions $g_\tau$ on the 2-sphere and
  $\phi_\tau\in H^1(W, \cv)$ such that
  $\{\phi_1, \phi_2\}=\sum_{\tau}\dbar\left( g_\tau\phi_\tau\right).$
  \end{lemma}
  \bproof
  We repeat part of a previous computation as follows.
  \begin{eqnarray*}
  & &\{ \lambda_1^k\lambda_2^{2-k}V_i^a\Oomega_j^b,
      \lambda_1^{k'}\lambda_2^{2-k'}V_{i'}^{m+1}\Oomega_{j'}^{b'}\}
  \\
  &=& \frac{\mu^{k+k'}}{(1+|\mu|^2)^2}
   \left(
  \partial_i^a
  \osigma_{j'}^{b'}\wedge L_{\partial_{i'}^{m+1}}(\osigma_j^b)
  +
  \partial_{i'}^{m+1}
  \osigma_j^b \wedge L_{\partial_{i}^{a}} (\osigma_{j'}^{b'})
  \right)
  \\
  &=& \frac{\mu^{k+k'}}{(1+|\mu|^2)^2}
   \left(
  \partial_{i'}^{m+1}
  \osigma_j^b\wedge L_{\partial_{i}^{a}}(\osigma_{j'}^{b'})
   \right)
  = \frac{\mu^{k+k'}}{(1+|\mu|^2)^2}
  \partial_{i'}^{m+1}
  \osigma_j^b\wedge \delta_{ab'}\epsilon_{ij'}\frac{d{\mubar}}{1+|\mu|^2}
  \\
  &=&
  -\frac{\mu^{k'}}{(1+|\mu|^2)^2}
  (\delta_{ab'}\epsilon_{ij'}\frac{\mu^k}{1+|\mu|^2}\partial_{i'}^{m+1}
  \osigma_j^b)
  =
  -\dbar(f\lambda_1^k\lambda_2^{2-k}\delta_{ab'}\epsilon_{ij'}
  V_{i'}^{m+1}\Oomega_j^b),
  \end{eqnarray*}
  where $f$ is one of the functions $f_1$, $f_2$ and $f_3$ depending on
  the number $k'$. Since
  $\lambda_1^k\lambda_2^{2-k}V_{i'}^{m+1}\Oomega_j^b$
  is in $H^1(W, \cv)$, the proof of this lemma is
  completed.
  \eproof

  \begin{lemma}For any $\phi_1$ in $ker\delta_1$ and
  $\phi_2=\lambda_1^{k'}\lambda_2^{2-k'}V_{i'}^{m+1}\Oomega_{j'}^{m+1}$,
  $\{\phi_1, \phi_2\}=0.$
  \end{lemma}
  \bproof We have six cases to consider.
  \begin{eqnarray*}
  & &\{ \lambda_1^k\lambda_2^{2-k}(V_1^a\Oomega_2^b+V_1^b\Oomega_2^a),
      \lambda_1^{k'}\lambda_2^{2-k'}V_{i'}^{m+1}\Oomega_1^{m+1}\}
  \\
  &=&
  \frac{\mu^{k+k'}\partial_{i'}^{m+1}}{(1+|\mu|^2)^2}
  \left(
  \osigma_2^b \wedge L_{\partial_1^a}\osigma_1^{m+1}
  +   \osigma_2^a \wedge L_{\partial_1^b}\osigma_1^{m+1}
    \right)
  \\
  &=&
  2{\mubar}\frac{\mu^{k+k'}}{(1+|\mu|^2)^2}\partial_{i'}^{m+1}
    \left(
    \osigma_2^b\wedge\osigma_2^a+\osigma_2^a\wedge\osigma_2^b
  \right)=0.
  \end{eqnarray*}
  Similarly,
  \[
  \{ \lambda_1^k\lambda_2^{2-k}(V_2^a\Oomega_1^b+V_2^b\Oomega_1^a),
      \lambda_1^{k'}\lambda_2^{2-k'}V_{i'}^{m+1}\Oomega_1^{m+1}\}
  =
  \frac{-2{\mubar}\mu^{k+k'}\partial_{i'}^{m+1}
    \left(
    \osigma_1^b\wedge\osigma_1^a+\osigma_1^a\wedge\osigma_1^b
  \right)}{(1+|\mu|^2)^2}.
  \]
  It is equal to zero. Next,
  \begin{eqnarray*}
  & &\{ \lambda_1^k\lambda_2^{2-k}(V_1^a\Oomega_1^b-V_2^b\Oomega_2^a),
      \lambda_1^{k'}\lambda_2^{2-k'}V_{i'}^{m+1}\Oomega_1^{m+1}\}
  \\
  &=&
  \frac{\mu^{k+k'}}{(1+|\mu|^2)^2}\partial_{i'}^{m+1}
    \left(
  \osigma_1^b \wedge L_{\partial_1^a}\osigma_1^{m+1}
  -  \osigma_2^a \wedge L_{\partial_2^b}\osigma_1^{m+1}
    \right)
  \\
  &=&
  2{\mubar}\frac{\mu^{k+k'}}{(1+|\mu|^2)^2}\partial_{i'}^{m+1}
    \left(
    \osigma_1^b\wedge\osigma_2^a+\osigma_2^a\wedge\osigma_1^b
  \right)=0.
  \end{eqnarray*}
  Similarly,
  \[
  \{ \lambda_1^k\lambda_2^{2-k}(V_1^a\Oomega_1^b-V_2^b\Oomega_2^a),
      \lambda_1^{k'}\lambda_2^{2-k'}V_{i'}^{m+1}\Oomega_2^{m+1}\}
  =\frac{2\mu^{k+k'}\partial_{i'}^{m+1}
    \left(
    \osigma_1^b\wedge\osigma_2^a+\osigma_2^a\wedge\osigma_1^b
  \right)}{(1+|\mu|^2)^2}.
  \]
  It is equal to zero. Finally,
  \begin{eqnarray*}
  & &\{ \lambda_1^k\lambda_2^{2-k}(V_1^a\Oomega_2^b+V_1^b\Oomega_2^a),
      \lambda_1^{k'}\lambda_2^{2-k'}V_{i'}^{m+1}\Oomega_2^{m+1}\}
  \\
  &=&
  \frac{\mu^{k+k'}}{(1+|\mu|^2)^2}\partial_{i'}^{m+1}
    \left(
  \osigma_2^b \wedge L_{\partial_1^a}\osigma_2^{m+1}
  +   \osigma_2^a \wedge L_{\partial_1^b}\osigma_2^{m+1}
    \right)
  \\
  &=&
  2\frac{\mu^{k+k'}}{(1+|\mu|^2)^2}\partial_{i'}^{m+1}
    \left(
    \osigma_2^b\wedge\osigma_2^a+\osigma_2^a\wedge\osigma_2^b
  \right)=0.
  \end{eqnarray*}
  Similarly, the following term is equal to zero because
  \[
  \{ \lambda_1^k\lambda_2^{2-k}(V_2^a\Oomega_1^b+V_2^b\Oomega_1^a),
      \lambda_1^{k'}\lambda_2^{2-k'}V_{i'}^{m+1}\Oomega_2^{m+1}\}
  =\frac{-2\mu^{k+k'}\partial_{i'}^{m+1}
    \left(
    \osigma_1^b\wedge\osigma_1^a+\osigma_1^a\wedge\osigma_1^b
  \right)}{(1+|\mu|^2)^2}.
  \]
  \eproof

  All lemmas in this section together prove
  that when $\phi_1$ and $\phi_2$ are
  in $\ce$, then there exist
  functions $g_\tau\in \Gamma^0$ such that
  $
  \{\phi_1, \phi_2\}=\sum \dbar \left(g_\tau\Up_\tau \right).
  $
  Suppose $\phi_1, \phi_2$ are in $\ce$
  and $h_1, h_2$ are in $\Gamma^0$. Since
  $\Up_\tau (g)=0$ if  $g\in\Gamma^0$,
  \[
  \{h_1\phi_1, h_2\phi_2\}=h_1h_2\{\phi_1, \phi_2\}
  =\sum_\tau h_1h_2 \dbar \left(g_\tau\Up_\tau \right)
  =\sum_\tau \left(h_1h_2 \dbar g_\tau\right)\wedge\Up_\tau .
  \]
  Since $h_1h_2 \dbar g_\tau$ is a $\dbar$-closed (0,1)-form on
  the Riemann sphere, it is $\dbar$-exact.
  Our proof of Proposition \ref{expansion} is now completed.

  \subsection{Convergence}

  Choose a Hermitian inner product on $\ce$
  such that the elements $\Up_\tau$ form a Hermitian basis.
  On the Sobolev spaces $L_k^2(S^2)$ and
  $L_k^2(\Gamma^{(0,1)}(S^2))$, we have the
  usual  quadratic norm $\|\cdot\|_k$. Define a norm on
  $\| \cdot \|_k$ on $\Gamma^0\otimes {\ce}$
  by $\|\sum_\tau f_\tau\Up_\tau\|^2_k =\sum_\tau\|f_\tau\|^2_k$.
  Similarly, define a norm
  $\|\cdot \|_k$ on $\Gamma^{(0,1)}(S^2)
  \otimes {\ce}$ by
  $\|\sum_\tau\gamma_\tau \wedge \Up_\tau\|^2_k
  =\sum_\tau\|\gamma_\tau\|^2_k$.

  \begin{lemma}\label{Schauder} There exists a constant $C_1$ such that
  if $f$ is a smooth function on $S^2$ with $\int_{S^2} f =0$, then
  $\|f\|_k \leq C_1 \|\dbar f\|_{k-1}.$
  \end{lemma}
  \bproof  The condition on $f$ implies that it is
  orthogonal to $\ker\dbar^*\dbar$.
  By Schauder estimate \cite[Appendix H,
  Theorem 27]{Besse},
   there exists a constant
  $c_1$ such that for all such $f$,
  $
  \|f\|_k \leq c_1\|\dbar^* \dbar f\|_{k-2}.
  $

  Consider the first order elliptic operator
  $\dbar + \dbar^* : \Gamma^{(0,1)} \to  \Gamma^{(0,0)} \oplus
  \Gamma^{(0,2)}$. Now $\dbar^* \dbar f= (\dbar + \dbar^*)\dbar f$.
  Since $\dbar f$ is orthogonal to $\ker (\dbar + \dbar^*)$,
  again by Schauder estimate
  there exists a constant
  $c_2$ such that for all such $f$,
  $
  \|\dbar^* \dbar f\|_{k-2}\leq c_2 \|\dbar f\|_{k-1}.
  $
  The Lemma follows.
  \eproof

  For the following lemma, see also \cite[5.118, 5.119]{Kodaira}.

  \begin{lemma}{\rm \label{estimate}}
  There exist constants $C_2$ and $C_3$ such that
  for all $\psi$ and $\phi$ in
  $\Gamma^0\otimes \ce$,
  \begin{eqnarray}
  \|\phi\|_k &\leq& C_2\|\dbar \phi\|_{k-1}; \label{first}\\
  \|\{\psi,\phi\}\|_k &\leq& C_3 \|\psi\|_{k+1}\| \phi\|_{k+1} \label{second}.
  \end{eqnarray}
  \end{lemma}
  \bproof
  Since $\{\Up_\tau: 1\leq \tau\leq D\}$ is a Hermitian basis and elements are holomorphic,
   (\ref{first}) follows from Lemma \ref{Schauder}.

  To prove (\ref{second}),  we
  assume that $\psi=g_1A_1$, $\phi=g_2A_2$ where $A_1$ and
  $A_2$ are one of the elements in the Hermitian basis $\{\Up_\tau: 1\leq
  \tau \leq D\}$.
  By Lemma \ref{expansion}, there exists $f_\tau$ such that
  $\{g_1A_1,g_2A_2\} = g_1g_2 \sum_\tau \dbar f_\tau \wedge \Up_\tau.$
  By definitions,
  $\|\{g_1A_1,g_2A_2\}\|_k^2 \leq
  \sum_\tau \|g_1g_2 \dbar f_\tau\|_k^2.$
  Define
  \[
  c_2(A_1, A_2)=\sqrt{\max_{\ell\leq k, \tau}\sup |\nabla^\ell\dbar f_\tau |^2}.
  \]
   This constant depends on  $A_1$ and $A_2$ because
  $f_\tau$ does. Define
  \[
  c_2=\max \{c_2(A_1, A_2): A_1, A_2\in\{\Up_\tau: 1\leq \tau\leq D\}\}.
  \]
  Then
  $\|g_1g_2 \dbar f_\tau\|_k \leq c_2 \|g_1g_2\|_k.$
  It is known \cite[page 73]{Kur2} that there exists a constant $c_3$ such that
  $\|g_1g_2\|_k \leq c_3\|g_1\|_k\|g_2\|_k$
  if $k \geq \dim S^2+1=3$.
  Combining all these inequalities, we have
  \begin{eqnarray*}
  \|\{g_1A_1,g_2A_2\}\|_k &\leq& D c_2 \|g_1g_2\|_k \leq
  D c_2c_3 \|g_1\|_k \|g_2\|_k\\
  &=& D c_2c_3
  \|g_1A_1\|_k\|g_2A_2\|_k \leq D c_2c_3
  \|g_1A_1\|_{k+1}\|g_2A_2\|_{k +1}.
  \end{eqnarray*}

  In general, $\psi=\sum_\tau\psi_\tau \Up_\tau$, $\phi=\sum_\tau\phi_\tau\Up_\tau$.
  For $k \geq \dim S^2+1=3$, we obtain:
  \begin{eqnarray*}
  & &\|\{\psi, \phi\}\|_k^2
    =\|\{\sum \psi_\tau\Up_\tau , \sum \phi_\rho\Up_\rho\}\|^2_k
  \leq \sum_{\tau, \rho} \|\{\psi_\tau\Up_\tau , \phi_\rho\Up_\rho\}\|^2_k\\
   &\leq& D^2c^2_2c_3^2\sum_{\tau, \rho}
  \|\psi_\tau\|^2_{k+1}\|\phi_\rho\|^2_{k+1}
  =D^2c^2_2c_3^2
  (\sum_{\tau}
  \|\psi_\tau\|^2_{k+1} )
  ( \sum_{\rho}\|\phi_\rho\|^2_{k+1})\\
  &=&D^2c^2_2c_3^2 \|\sum\psi_\tau\Up_\tau\|^2_{k+1}\cdot
  \|\sum\phi_\rho\Up_\rho\|^2_{k+1}
  =  D^2c^2_2c_3^2\|\psi\|^2_{k+1}\|\phi\|^2_{k+1}
  \end{eqnarray*}
  The proof is completed. \eproof

  \begin{lemma}\label{convergence}
  Let $\Phi = \sum_n\Phi_n t^n$. Suppose that the
  coefficients $\Phi_k$ are vector-valued $(0,1)$-forms with the property that
  \begin{equation}\label{recursive}
  \dbar\Phi_{n+1} = -1/2 \sum_{i=1}^n \{\Phi_i, \Phi_{n+1-i}\}
  \end{equation}
   and suppose they satisfy
  the inequalities {\rm (\ref{first})} and
  {\rm (\ref{second})}. Then  for small enough $t$,
  $\Phi$ is convergent.
  \end{lemma}
  \bproof Take a power
  series $A(t) = \sum_{n=1}^\infty a_nt^n = \sum k\frac{(ct)^n}{n^2}$.
  It has nonzero radius of convergence. As seen in the proof of
  \cite[Equation 5.116]{Kodaira},
  \[
  \frac{1}{16k} \sum_{i=1}^n a_ia_{n+1-i} \leq a_{n+1}.
  \]
  Suppose that $\|\Phi_n\|^2_k \leq a_n $ up to some $n$. By the last lemma,
  \[
  \|\Phi_{n+1}\|^2_k \leq  C^2_2\|\dbar\Phi_{n+1}\|^2_{k-1}
    \leq \frac{C^2_2}{4}\sum_{i=1}^n \|\{\Phi_i, \Phi_{n+1-i}\}\|^2_{k-1}
   \leq  \frac{C^2_3C^2_2}{4}\sum_{i=1}^n \|\Phi_{n+1-i}\|^2_k\|\Phi_i\|^2_k.
  \]
  If we choose $\frac{1}{16k} = \frac{C^2_3C^2_2}{4}$ and $c$ so that $\|\Phi_1\|^2\leq
  a_1=kc$, we obtain:
  \[
  \|\Phi_{n+1}\|^2_k \leq \frac{1}{16k} \sum_{i=1}^n \|\Phi_{n+1-i}\|^2_k
  \|\Phi_i\|^2_k
  \leq \frac{1}{16k} \sum_{i=1}^n a_{n+1-i}a_{i}
  \leq a_{n+1}
  \]
  So the Lemma follows by induction. The convergency of $A(t)$ is
  proved.
  \eproof

  For any element $\Phi_1$ in $H^1(W, \cd_W)\subset H^1(W, \Theta_W)$,
  Lemma \ref{expansion} inductively determines solutions for the recursive formula
  $\dbar\Phi(t)+\frac12\{\Phi(t), \Phi(t)\}=0.$
  Lemma \ref{convergence} shows that there exists $\epsilon >0$
  such that the power series
  $\Phi(t)=\Phi_1 t+\Phi_2 t^2+ \dots$
  converges when $|t|<\epsilon$.
  It follows from Kodaira-Spencer theory that $\Phi_1$ is an infinitesimal
  deformation. Proof of Theorem \ref{theorem2} is completed.

  \subsection{Hypercomplex Deformations}\label{hyper D}
  The integrability of $\Phi_1$ amounts to
  a deformation of quaternionic structures \cite{PP1}.
  We complete our discussion on deformation of hypercomplex structures
  by showing that every quaternionic deformation of the given hypercomplex structure
  on $X$ is a hypercomplex deformation. This is Theorem \ref{theorem3}.

  The underlying smooth structure of
  the twistor spaces of the deformed quaternionic structures on
  $X$ remains to be $X\times S^2$. The underlying
  smooth structure of the twistor spaces
  of the deformed quaternionic structures on $T^{4m}$ is $T^{4m}\times S^2$.
  The natural projections $p_W$ and $p_Z$ from these spaces
  onto $S^2$ satisfy the identity $p_W=p_Z\circ \Psi$.
  When the twistor spaces are given the un-perturbed hypercomplex structures,
  the maps $p_W$ and $p_Z$ are holomorphic. We complete the proof of
  Theorem \ref{theorem3} by showing that for all local deformations found in
  the previous paragraph, the maps $p_W$ and $p_Z$
  are holomorphic. They are denoted by $p$ in the previous sections.
  As a result of general twistor theory \cite{PP1} \cite{PP2},
  the deformed quaternionic manifolds are in fact hypercomplex manifolds.

  To verify our claim on the open set $U_2$,
  we consider any distribution determined by the convergent power
  series $\Phi(t)=\sum_n \Phi_nt^n$ where $\Phi_1$ is in $H^1(W, \cd_W)\subset\ce$,
  and $\Phi_n$ is in $\Gamma^0\otimes\ce$ for all $n\geq 2$.
  On the open set $U_2$ of the twistor space $W$, the space of $(0,1)$-vectors is
  spanned by $\{{\overline V}_i^\alpha, \frac{\partial}{\partial \mubar}\}$.
  Then $(0,1)$-vectors for the deformation family are spanned by
  \begin{equation}
  {\overline V}_k^{\gamma}(t)
  :={\overline V}_k^\gamma
  +\Phi(t)({\overline V}_k^\gamma)
  ={\overline V}_k^\gamma+\sum_nt^n\Phi_n({\overline V}_k^\gamma)
  \end{equation}
  and $\frac{\partial}{\partial {\overline\mu}}$. To prove that $p_W$ is holomorphic for this
  family of complex structures, we verify that its differential
   sends $(0,1)$-vectors to
  $(0,1)$-vectors. Since $dp_W(\frac{\partial}{\partial {\overline\mu}})=
  \frac{\partial}{\partial {\overline\mu}}$, our only concerns are on
  $dp_W{\overline V}_k^{\gamma}(t)$
  which is equal to
  \begin{equation}
  dp_W({\overline V}_k^\gamma)+\sum t^ndp_W(\Phi_n({\overline V}_k^\gamma)).
  \end{equation}

  By construction, $\Phi_n$ is in the space $\ce$ (\ref{spaceE}). Therefore,
  $\Phi_n({\overline V}_k^\gamma)$ is a linear combinations of
  $V_\ell^\rho$ for some $\ell$ and $\rho$.
  Therefore, it suffices to prove that $dp_W({\overline V}_k^\gamma)=0$ for
  all $k$ and $\gamma$. Indeed, when $\gamma=m+1$, ${\overline V}_k^\gamma$
  is given by coordinate vector fields $\frac{\partial}{\partial
  {\overline w}_i^{m+1}}$.
  As these vector fields are vertical with respect to the project $p_W$,
  they are in the kernel of the differential of $p_W$.
  Strictly speaking, when $1\leq \gamma\leq m$, the element
  ${\overline V}_k^\gamma$ is contained in $H^0(W, \Psi^*{\cd})$.
  They could be interpreted as tangent vectors on $W$ only after a
  lifting process as described in  our coboundary map computation
  in the proof of Theorem \ref{theorem2}.
  In particular, as seen in (\ref{lifting})
  the lifting of $\frac{\partial}{\partial x_j}$ is
  $X_j$ and $\frac{\partial}{\partial y_j}$ is $Y_j$.
  Taking this lifting into account, we have the identity
  \[
  dp_W({\overline V}_k^\gamma)
  =dp_Z\circ d\Psi ({\overline V}_k^\gamma)
  =dp_Z({\overline V}_k^\gamma)
  \]
  where ${\overline V}_k^\gamma$ on the right hand side of the equation is now
  interpreted as a vector field on $Z$.
  On $Z$, ${\overline V}_k^\gamma$ is spanned by
  $\frac{\partial}{\partial {\overline w}_j^\gamma}$. As it is obviously
  in the kernel of $dp_Z$, we complete the proof that the map $p_W$ is holomorphic.

  \subsection{Quaternionic deformations}
  Computation in the last section suggests that we have to
  take a closer look at the parameter count for hypercomplex
  and quaternionic deformations and their relations.
  The relation could be manifested by the details in the following
  result.
  \begin{lemma}  The coboundary maps $\delta_0$ and $\delta_1$ from the
  zero-th and first cohomology respectively in the induced cohomology sequence of
  \begin{equation}
  0 \to \cd_W \to \Theta_W   \to p^*\co(2)\to 0
  \end{equation}
  are injective.
  \end{lemma}
  \bproof
  Due to Lemma \ref{cohomology},
  the following twisted k-forms form a basis for  $H^k(W, p^*\co (2))$.
  \begin{equation}
  \lambda_1^{2+k-l}\lambda_2^l
    \Oomega_1^{\alpha_1}\wedge\cdots\wedge\Oomega_1^{\alpha_i}
        \wedge\Oomega_2^{\beta_1}\wedge\cdots\wedge\Oomega_2^{\beta_j}
  \end{equation}
  where $0\leq l\leq 2+k$ and $1\leq \alpha_i, \beta_j\leq m+1$.

  Using the standard identification between quadratic polynomials and global holomorphic
  vector fields on $\cp^1$, i.e. $(a\lambda_1^2+b\lambda_1\lambda_2+c\lambda_2^2)
  \mapsto (a\mu^2+b\mu+c)\frac{\partial}{\partial\mu}$, we
  consider the above sections as holomorphic 1-forms on the twistor space with
  values in the tangent bundle of $\cp^1$. Then the above sections are identified to
  \begin{equation}
  s=\mu^{2+k-l}\frac{\partial}{\partial\mu}
    \osigma_1^{\alpha_1}\wedge\cdots\wedge\osigma_1^{\alpha_i}
        \wedge\osigma_2^{\beta_1}\wedge\cdots\wedge\osigma_2^{\beta_j}.
  \end{equation}
  To find the image of $s$ through the coboundary map $\delta_k$,
  we lift the 1-forms by pull-backs and the vector field is lifted to
  $W$. It defines a smooth lifting $\hat{s}$.
  We use the Chern connection $\nabla$
  on the holomorphic tangent bundle on $W$ to define the $\dbar$-operator
  $\dbar^\nabla$. Then $\delta_k s$ is represented uniquely by the cohomology
  class of $\dbar^\nabla \hat{s}$. Since the function $\mu^{\ell+2-l}$  and the forms
  $\osigma_1^{\alpha}$ and $\osigma_2^{\beta}$ are holomorphic,
  \begin{equation}
  \dbar^\nabla \hat{s}=
  \mu^{2+k-l}\left(\dbar^\nabla\frac{\partial}{\partial\mu}\right)
    \wedge\osigma_1^{\alpha_1}\wedge\cdots\wedge\osigma_1^{\alpha_i}
        \wedge\osigma_2^{\beta_1}\wedge\cdots\wedge\osigma_2^{\beta_j}.
  \end{equation}
  Due to an observation of Gauduchon \cite{Gauduchon},
  for any (0,1)-vector $X$, and (1,0)-vector $Y$,
  $\dbar^\nabla_X Y=[X, Y]^{1,0}$.
  Using (\ref{32}) and (\ref{34}), we find that
  \[
  [\dbar_1^\alpha,\frac{\partial}{\partial\mu}]^{1,0} = -\frac{\partial_2^\alpha}{1+|\mu|^2},
  \qquad
  [\dbar_2^\alpha,\frac{\partial}{\partial\mu}]^{1,0} = \frac{\partial_1^\alpha}{1+|\mu|^2},
  \]
  It follows that $\dbar^\nabla\left(\frac{\partial}{\partial\mu}\right)$
  is equal to
  \begin{eqnarray*}
  &\ & -\frac{d\wbar^\alpha_1\otimes\partial_2^\alpha}{1+|\mu |^2}     +\frac{d\wbar^\alpha_2\otimes\partial_1^\alpha}{1+|\mu |^2}
   +\frac{w_1^\alpha-\mu\wbar^\alpha_2}{(1+|\mu |^2)^2}d\mubar_1\otimes\partial_1^\alpha
    +\frac{w_2^\alpha+\mu\wbar^\alpha_1}{(1+|\mu |^2)^2}d\mubar_1\otimes\partial_2^\alpha
   \\
  &=& -\frac{d\wbar^\alpha_1\otimes\partial_2^\alpha}{1+|\mu |^2}     +\frac{d\wbar^\alpha_2\otimes\partial_1^\alpha}{1+|\mu |^2}
   -\frac{\zbar^\alpha_2}{(1+|\mu |^2)^2}d\mubar_1\otimes\partial_1^\alpha
    +\frac{\zbar^\alpha_1}{(1+|\mu |^2)^2}d\mubar_1\otimes\partial_2^\alpha
   \\
  &=&\sum_{\alpha}\left(\osigma^\alpha_2\otimes\partial_1^\alpha-\osigma^\alpha_1\otimes\partial_2^\alpha
  \right).
  \end{eqnarray*}
  Therefore, $\dbar^\nabla_k \hat{s}$ is equal to
  \begin{eqnarray*}
   & &
  \mu^{2+k-l}
  \sum_{\alpha}
  \left(\osigma^\alpha_2\otimes\partial_1^\alpha-\osigma^\alpha_1\otimes\partial_2^\alpha
  \right)
    \wedge\osigma_1^{a_1}\wedge\cdots\wedge\osigma_1^{a_i}
        \wedge\osigma_2^{b_1}\wedge\cdots\wedge\osigma_2^{b_j}
   \\
  &=& \lambda_1^{2+k-l}\lambda_2^l
  \sum_{\alpha}\left(\Oomega^\alpha_2\otimes\frac{1}{\lambda_2}\partial_1^\alpha-
    \Oomega^\alpha_1\otimes\frac{1}{\lambda_2}\partial_2^\alpha
  \right)
    \wedge\Oomega_1^{a_1}\wedge\cdots\wedge\Oomega_1^{a_i}
        \wedge\Oomega_2^{b_1}\wedge\cdots\wedge\Oomega_2^{b_j}
  \\
  &=& \lambda_1^{2+k-l}\lambda_2^l
  \sum_{\alpha}\left(\Oomega^\alpha_2\otimes V^\alpha_1-
    \Oomega^\alpha_1\otimes V^\alpha_2
  \right)
    \wedge\Oomega_1^{a_1}\wedge\cdots\wedge\Oomega_1^{a_i}
        \wedge\Oomega_2^{b_1}\wedge\cdots\wedge\Oomega_2^{b_j}.
  \end{eqnarray*}
  This map is injective when $k=0,1$.
  In particular,
  \begin{equation}
  \delta_0(\lambda^{2-l}_1\lambda^l_2)=
  \lambda_1^{2-l}\lambda_2^l
  \sum_{\alpha}\left(\Oomega^\alpha_2\otimes V^\alpha_1-
    \Oomega^\alpha_1\otimes V^\alpha_2
  \right)
  \end{equation}
  The proof is now completed.
  \eproof

  As a consequence of the last lemma, the following sequence
  is exact
  \begin{equation}\label{surjection}
  0\to H^0(W, p^*\co(2))\to H^1(W, \cd_W)\to H^1(W, \Theta_W)\to 0.
  \end{equation}
  It confirms, but does not prove our result
  in the last section, that every quaternionic deformation is a
  hypercomplex deformation. The new information in the above
  exact sequence is, that for every hypercomplex deformation parameter,
  there is a three-dimensional hypercomplex deformation within one
  quaternionic class. This parameter space is contributed by
  $H^0(W, p^*\co(2))\cong S^2\bC^2$. Since $H^1(W, \Theta_W)$ is
  the virtual parameter space for complex structures on the twistor
  space $W$, its real subspace with respect to the real structure $\tau$
  is the parameter space for quaternionic structure on $X$, it completes
  the proof of Theorem \ref{theorem4}.

  \subsection{Deformations of Torus}
  During the course of our computation, significant amount of information on deformation
  of the hypercomplex structure on the torus $T^{4m}$ as a quotient of
  the quaternion module $\bH^m$ is revealed.
  By Lemma \ref{torus moduli}, the virtual parameter space for hypercomplex deformation
  on the torus $T^{4m}$ is the real part of the complex vector space $\lie{t}^{1,0}\otimes
   \lie{t}^{*(0,1)}\otimes S^{2}\bC^{2}$. Its real dimension is equal to complex dimension;
  i.e. $(2m)^2\times 3=12m^2$. To prove integrability, we apply the method
  in Section \ref{integrable} with the space $\ce$ replaced by $H^1(Z, \cd_Z)$. In this case,
  we do not need the work in
  Section \ref{hyper D} to conclude that every deformation generated by $H^1(Z, \cd_Z)$ is
  hypercomplex as the power series generated in Section \ref{integrable} has values in
  $H^1(Z, \cd_Z)$.
  Finally, the exact sequence (\ref{surjection}) has its counter part on $Z$.
  Therefore, the parameter space for quaternionic deformation on $T^{4m}$ is equal to
  $12m^2-3$.

  In Section \ref{defining}, we construct a hypercomplex structure on $\bR^{4m}$
  by left multiplication of unit quaternions $i$, $j$ and $k$ on $\bH^m$.  As
  seen in Section \ref{algebra}, the torus $T^{4m}$ is the quotient of $\bR^{4m}$
  with respect to the lattice group $\phi(\Gamma)\cong (\bZ^{4m}, +)$.  The right
  multiplication by a generic element in $\GL(4m, \bR)$ changes the hypercomplex
  structure on $T^{4m}$ by choosing different identification from
  $\bR^{4m}$ to $\bH^{m}$. The isotropy subgroup with respect to the quaternion
  basis $\{Q_a, 1\leq a\leq m\}$ is $\GL(m, \bH)$. Therefore, the (coarse topological) moduli space of
  hypercomplex structures is $Aut(\bZ^{4m})\backslash \GL(4m, \bR)/\GL(m, \bH)$ where $Aut(\bZ^{4m})$ is the
  discrete subgroup of $GL(4m, \bR)$ consisting of
  automorphisms of the lattice. Furthermore, left
  multiplication of $\SP(1)$ changes the hypercomplex structure but keeps the
  quaternionic structure.
  Therefore, the (coarse topological) moduli space of quaternionic
  structures
  is $Aut(\bZ^{4m})\backslash \GL(4m, \bR)/\GL(m, \bH)\SP(1)$.
  Our computation on the dimension of
  moduli proves that these spaces are the entire connected component of the
  moduli space of hypercomplex deformations and quaternionic
  deformations respectively.

 \section{Moduli of Invariant Hypercomplex Structures}

By invariant hypercomplex structure we mean a triple of
left-invariant complex structures on $X$ with the usual relations.
Such triple is determined by its values at the identity. The first
observation in this section is due to the power series
calculations in Section 4.1.
\begin{proposition}\label{proposition2}
Any small deformation of the hypercomplex structure on $X$
consists of invariant structures.
\end{proposition}

\bproof From the calculations in the proof of Proposition 1,
Section 4.1  we can conclude that in the series $\Phi(t)= \sum_n
\Phi_n t^n$ all terms $\Phi_n$ belong to $\Gamma_0\otimes \ce$. In
particular for every fixed $t_0$, $\Phi(t_0)$ is in the same
space. Then substituting $(\lambda_1,\lambda_2)$ with $(1,0),
(0,1)$ and $(i,1)$ we obtain elements in $\ce$. They determine
invariant complex structures by identifying $\ce$ with a subspace
of the space of sections in $T^{(1,0)X}\otimes T^{*(0,1)}$. Each
section is identified with the image of the bijective linear map $
Id + \Phi :T_I^{(0,1)}X\rightarrow T_J^{(0,1)}X$ and the deformed
structure $J$ is defined from here. Since the elements of $\ce$
are given by invariant sections, the deformed structures are
invariant. Then the three invariant complex structures give rise
to an invariant hypercomplex structure. \eproof

  Further in this section we investigate the deformations of the hypercomplex structures on
   $(H_{4m+1}\times{\bf R}^3)/\Gamma$ arising from the deformations of the lattice
  $\Gamma$. As we will see this space differs from the space of
  invariant hypercomplex structures.
  %The space $(H_{4m+1}\times{\bf R}^3)/\Gamma$ is a hypercomplex
  %manifold because the hypercomplex structure on $H_{4m+1}\times{\bf
  %R}^3$ is left invariant.
%  We first calculate the exponential map
 % $exp: {\lie{h}}_{4m+1}\oplus {\bf R}^3 \rightarrow H_{4m+1}\times {\bf
 % R}^3$. Consider again the left invariant vector fields $X_i =
 % \frac{\partial}{\partial x_i} +2y_i\frac{\partial}{\partial z},
 % Y_i =\frac{\partial}{\partial y_i} -2x_i\frac{\partial}{\partial
 % z}, Z= \frac{\partial}{\partial z}, E_j = \frac{\partial}{\partial
 % t_j}$.  Let $V=\Sigma a_iX_i+b_iY_i+cZ+e_jE_j$. Then
 % $\phi(t)=exp(tV)$ is defined by:
 % $$
 % \frac{d\phi}{dt}(t)=V(\phi(t)) \hspace{.1in} \phi(0)=0
 % $$
 % From here we obtain
 % $
 % exp(V)=\phi(1)=(a_i,b_i,c,e_j).
 % $
 % It follows that the exponential map is a
 % diffeomorphism from ${\bf R}^{4n+4} = {\bf H}^{n+1}$ to
 % $H_{4m+1}\times \bf{R}^3$ which is represented by the identity in
 % the coordinates with respect to the given basis and the cannonical
 % one of ${\bf R}^{4m+4}$.

  From Theorem \ref{identification}  the left invariant hypercomplex structure
defined in Section 2.2  is equivalent to the canonical one on
${\bf H}^m$.

   Now if $\Gamma '$ is another lattice in $H_{4m+1}\times {\bf R}^3$ such that the
   quotient spaces $H_{4m+1}\times {\bf R}^3/\Gamma$ and
  $H_{4m+1}\times {\bf R}^3/\Gamma '$ are isomorphic as hypercomplex manifolds, then
  $\Gamma$ and $\Gamma '$ are isomorphic. According to (\cite{R},
  Theorem 2.11, Corollary 2) the isomorphism between the two cocompact subgroups
   is uniquely extended to an automorphism
  $\overline{\Upsilon}$ of $H_{4m+1}\times {\bf R}^3$. Then by the
  functorial property of the exponential map there is a unique
  Lie algebra automorphism $\Upsilon$ of $\lie{h}_{4m_+1}\oplus{\lie{t}}^3$
  defined by $\overline{\Upsilon}$. In particular the deformations arising from deformation
  of the lattice lie in the orbit of the automorphism group of $\lie{h}_{4m_+1}\oplus{\lie{t}}^3$
  containing the standard hypercomplex structure. Bellow we consider the form of
  $\Upsilon$.

  \begin{proposition}
  Let $\{Z,E_1,E_2,E_3,X_{2a-1},X_{2a},Y_{2a-1},Y_{2a} \}$ be an ordered basis for
  ${\lie{h}}_{4m+1}\times {\bf R}^3$. The automorphism group of $\lie{h}_{4m+1}\oplus \lie{t}^3$
  consists of elements $\Upsilon$
   leaving the center invariant of the form:
$$
  \Upsilon = \left( \begin{array}{cc}
               A & B \\
               0 & C \\
  \end{array}
  \right)
  $$
  where $A$ is in $End{\lie{c}}\oplus Hom (\lie{t}^3, \lie{c}\oplus\lie{t}^3)$,
  and $ B$ is in  $Hom(\lie{h}_{4m+1}\oplus \lie{t}^3, \lie{c}\oplus\lie{t}^3)$.
  Moreover, $C$ is a
  matrix preserving the symplectic form in ${\bf R}^{4m}$ up to a constant.
%  $$
 % \Upsilon = \left( \begin{array}{cccccccc}
 %              S_0 & S_1 &  S_2 & S_3 & F_{2a-1} & F_{2a} & G_{2a-1} & G_{2a} \\
 %              0 & 0 & 0 & 0 & P_{2a-1} & P_{2a} & Q_{2a-1} & Q_{2a} \\
 % \end{array}
 % \right)
 % $$
 % where $(S_0, S_1, S_2, S_3)$ is in $End{\lie{c}}\oplus Hom (\lie{t}^3, \lie{c}\oplus\lie{t}^3)$,
 % and $ (F_{2a-1}, F_{2a}, G_{2a-1}, G_{2a})$ is in  $Hom(\lie{h}_{4m+1}\oplus \lie{t}^3, \lie{c}\oplus\lie{t}^3)$.
 % Moreover, $(P_{2a-1}, P_{2a}, Q_{2a-1},Q_{2a})$ is a
 % matrix preserving a symplectic form in ${\bf R}^{4m}$ up to a constant.
  \end{proposition}
  \bproof The Lie
  brackets is
   $$[V,V'] = -2\omega(V,V')Z$$
  for the skew form $\omega$ with $Ker \omega =
  span\{Z,E_1,E_2,E_3\}$ and $\omega(V,V^{'})=
  2(-y_ix^{'}_i+x_iy^{'}_i)$ where $V=(z, e_j, x_i,y_i)$ and
  $V'=(z', e_j^{'}, x^{'}_i,y_i^{'})$ in the given basis. Since
  $\Upsilon$ preserves the center and the direct sum
  $\lie{h}_{4m+1}\oplus \lie{t}^3$, $\Upsilon(Z)=S_0 Z$ for some
  constant $S_0$. From the form of the brackets above we obtain
  $\omega(\Upsilon(V),\Upsilon(V'))=S_0\omega(V,V')$. \eproof

For further use denote the set of all $\Upsilon$ for which
   $\omega(\Upsilon(V),\Upsilon(V'))=S\omega(V,V')$ for a fixed $S$ as
$C_SSp(2m,{\bf R})$

  \begin{corollary}
  The dimension of the space of automorphisms in Proposition 2 is
  $13+18m+8m^2$
  \end{corollary}
  \bproof We calculate the dimensions of the three blocks. The
  $4\times 4$ block in the upper-left corner consisting of $S$'s has
  dimension 13. The upper-right block with $F_a$'s and $G_a$'s has
  dimension $4\times 4m = 16m$. Finally the set of symplectic
  transforms up to a constant $CSp(2m,{\bf R})$ has the same
  dimension as the group $Sp(2m,{\bf R})$ which is 2m(4m+1). When we
  sum up the dimensions the result follows.

\eproof

%{\bf Explanation why is $\dim
 % CSp(2m,{\bf R})=\dim Sp(2m,{\bf R})$?} \rm If $A$ is a matrix with
 % $\omega(A(V),A(U))=s\omega(V,U)$, then $\frac{A}{\sqrt{s} }$ is a
 % symplectic matrix for positive S. So $CSp(2m,{\bf R})$ is like a
 % left coset in GL and has the same dimension. When $s<0$, we have
 % another component of the automorphism group. But everything we
 % want is the dimension count - so its local.
%For  $s<0$ case we can first apply izomorphism A with $\omega(AX,AY)=-\omega(X,Y)$

  For any automorphism $\Upsilon$ of $h_{4m+1}\oplus\lie{t}^3$, we define
  hypercomplex structure ${\cal I}_{\Upsilon}$ by
  $\Upsilon^{-1}\circ {\cal I}\circ \Upsilon$ with the new basis $\Upsilon$.
  Then it will induce hypercomplex structure defined on
  $(H_{4m+1}\times {\bf R}^3)/{ \overline\Upsilon(\Gamma)}$ via
  factorisation. Now we have to consider which $\Upsilon$ give rise
  to equivalent structures on the factor-space.

  \begin{proposition}
  The hypercomplex spaces $(H_{4n+1}\times{\bf R}^3)/\Gamma$ and
  $(H_{4n+1}\times{\bf R}^3)/\overline{\Upsilon}(\Gamma)$ are
  equivalent if and only if:
\begin{eqnarray*}
  A & = & s I, s\neq 0 \\
  B &\in &CSp(2m, {\bf R})\cap GL(m,{\bf{H}})=CSp(m) \\
 C &=& (C_1  C_2 ... C_n), C_i =
   \left(
   \begin{array}{cccc}
   \alpha_{2i-1} &
  \alpha_{2i} & \beta_{2i-1} & \beta_{2i} \\
  -b_i & a_i & -d_i & c_i \\
  -c_i & d_i & a_i & -b_i \\
  -d_i & -c_i & b_i & a_i
  \end{array}
   \right)
  \end{eqnarray*}

%  \begin{eqnarray*}
%  (S_0,S_1,S_2,S_3) & = & s I, s\neq 0 \\
%  (P_{2a-1}, P_{2a}, Q_{2a-1},Q_{2a}) &\in &CSp(2m, {\bf R})\cap GL(m,{\bf{H}})=CSp(m) \\
%  (F_{2i-1}, F_{2i}, G_{2i-1}, G_{2i}) &=&
%   \left(
%   \begin{array}{cccc}
%   \alpha_{2i-1} &
%  \alpha_{2i} & \beta_{2i-1} & \beta_{2i} \\
%  -b_i & a_i & -d_i & c_i \\
%  -c_i & d_i & a_i & -b_i \\
%  -d_i & -c_i & b_i & a_i
%  \end{array}
%   \right)
%  \end{eqnarray*}
 % {\bf Could you please correct the indices and symbols in the
 % matrix above? Could you also explain how you get it? I am
 % confused. - The answer is at the end in the proof. I also changed
 % Sp to Csp in the statement above and the indices from $a$ to $i$}
  In particular the dimension of the group of hypercomplex
  automorphisms arising in this way is $1+9m+2m^2$.
  \end{proposition}
  \bproof
  Any hypecomplex automorphism of  $H_{4m+1}\times {\bf R}^3/\Gamma$
  gives rise to automorphism on $H_{4m+1}\times {\bf R}^3 = {\bf
  H}^{n+1}$ which interchange the actions of $\Gamma$ and $\Gamma'=
  \overline{\Upsilon}(\Gamma)$. But any hypercomplex automorphism of
  ${\bf H}^{n+1}$ is an affine transformation by a theorem of
  Ehresman. So $(A,{\bf v}): {\bf q} \rightarrow A{\bf q}+{\bf v}$
  with $A\in GL(n+1,{\bf H})$ is
   a general form of a hypercomplex automorphism. Then to verify the proposition we
  have to check that if
  $$(A,{\bf v})\circ \gamma = \gamma^{'}\circ (A,{\bf v})$$
  with $\gamma \in \Gamma$ and $\gamma^{'} =
  \overline{\Upsilon}(\gamma) \in \overline{\Upsilon}(\Gamma)$,
  then $\Upsilon$ has the proposed form.

  The group multiplication is expressed in terms of $\omega$ as  $$x
  *x^{'} = x+x^{'} + \omega(x,x^{'})z$$
   where $z=(1,,0,0,...,0)$ is the center generator
. Then the condition above becomes $A(\gamma * {\bf q}) + {\bf v}
=
  \gamma^{'} * (A{\bf q} + {\bf v})$. It follows that
  $$
  A\gamma + \omega(\gamma,{\bf q})Az= \gamma^{'} +
  \omega(\gamma^{'},A{\bf q}+{\bf v})z.
  $$
  Substituting ${\bf q}=0$, we have
  $$
  A\gamma = \gamma^{'} + \omega(\gamma^{'},{\bf v})z.
  $$
  Substitute the last formula back, we have
  $$
  \omega(\gamma^{'},A{\bf q})z = \omega(\gamma,{\bf q})Az
  \mbox{ and }
  A\gamma = \gamma^{'} + \omega(\gamma,A^{-1}{\bf v})Az.
  $$
  It follows that
  $$
  \overline{\Upsilon}(\gamma) = \gamma ' = A\gamma - \omega(\gamma,A^{-1}{\bf v})Az.
  $$
  As both $A$ and $\omega$ are linear, $\overline{\Upsilon}$ is
  linear. The above formula is applicable to any element in
  $H_{4m+1}\times {\bf R}^3$. It also shows that $(A, {\bf v})$
  uniquely determines $\overline{\Upsilon}$.

  Now using the identification of $\Upsilon$ and $\overline\Upsilon$
  via the exponential map we are going to obtain the form of  $\Upsilon$
  in the proposition. Note that A is simultaneuosly Lie-algebra
  automorphism described in Proposition 2 and an element of $GL(m+1,
  {\bf H})$.

 We first consider the case ${\bf v}=0$.
  Here we have to characterize the matrices $A$ as above.
   The
  group $GL(m+1,{\bf H})$ is identified with the group preserving
  the hypercomplex structure. The structures
  $I_1,I_2,I_3$ are formed by $4\times 4$-blocks along the diagonal
  of the following type:

  \begin{eqnarray*}
  J_1 &=& \left(
  \begin{array}{cccc}
   0 & -1 & 0 & 0 \\
  1 & 0 & 0 & 0 \\
  0 & 0 & 0 & -1 \\
  0 & 0 & 1 & 0
  \end{array}
  \right)
  \end{eqnarray*}

  and
  \begin{eqnarray*} J_2 &=& \left(
  \begin{array}{cccc}
   0 & 0 & -1 & 0 \\
  0 & 0 & 0 & 1 \\
  1 & 0 & 0 & 0 \\
  0 & -1 & 0 & 0
  \end{array}
  \right)
  \end{eqnarray*}
  and $J_3=J_1J_2$. Then when we divide the matrix of $\Upsilon$ in
  $4 \times 4$-blocks, then each block should represent a matrix
  which commute with $J_1$ and $J_2$ above.  The general form of a
  $4\times 4$ matrix which commutes with $J_1$ and $J_2$ is

  \begin{eqnarray*} B &=& \left(
  \begin{array}{cccc}
   a & b & c & d \\
  -b & a & -d & c \\
  -c & d & a & -b\\
  -d & -c & b & a
  \end{array}
  \right)
  \end{eqnarray*}
  for some numbers $a,b,c,d$. The diagonal blocks additionally have
  nonzero determinant.

  Now the upper-left corner of $A$ has first column with zero's
  except $S_0$, so it is proportional to the identity. The
  lower-right corner is $4m\times 4m$ consisting of
  $P_{2a-1},P_{2a},Q_{2a-1},Q_{2a}$ is an element of $C_{S_0}Sp(m)$. Finally
  the upper-right $4\times 4m$ corner consists of $4\times 4$
  matrices $(F_{2i-1}...G_{2i})$ for $i =1,...m$. Here $i$
  represents the number of the block. Then all
  matrices have the form of $B$ above.

  Now we have to consider the case ${\bf v} \neq 0$. In the formula
  $\gamma \rightarrow A\gamma -\omega(\gamma, A^{-1}{\bf v})Az$ we
  have that the second term is always proportional to $z$ i.e.
  $\omega(\gamma, A^{-1}{\bf v})Az = C z$. By varying $v$ we may
  choose any constant of proportion  $C$. Now when $\gamma$ is
  chosen to be any of the basic vectors the formula above changes an
  element of the first row of $\Upsilon$ because of this term. Thus
  the first row of $\Upsilon$ is arbitrary. With this we have the form
  of the matrix $(F_{2i-1},F_{2i},G_{2i-1},G_{2i})$ in the
  proposition since the center is the first element in the basis.
  From here the proposition follows.

  \eproof

  \begin{corollary}
  The dimensions of automorphisms in Proposition 3 is $1+9m+2m^2$
  \end{corollary}
  \bproof The calculation is similar to that in the previous
  Corollary. This time the dimension of upper-left corner is 1,
  because the matrix is proportional to identity. As efore the dimension of
  the set $CSp(m)$ is the same
%\bf WHY???? -  again, it is similar
%  to the answer of your first question. CSp(m) is a left coset of
%  the form $s.Sp(m)$ so has the same dimension\rm
  as the dimension of $Sp(m)=Sp(m,{\bf C})\cap U(2m) = m(2m+1)$.
  Finally the block with $F_a$'s and $G_a$'s has dimension $2\times
  4m = 8m$. This gives the corollary. \eproof

  From the above  we obtain that the space of effective parameters
  for the deformations of the hypercomplex structures on $X$ arising
  from deformation of the lattice is $\frac{G}{H}$ where $G$ is the
  group of all $\Upsilon$'s in Proposition 2 while $H$ is the group
  described in the Proposition 3. The dimension of this space is $12
  + 9m + 6m^2$.
%  \begin{remark}
 % In the proof we didn't identify the structure arising from
 % automorphisms of the lattice, which is a discrete subgroup of $G$.
 % So we have a deformation space instead of moduli space.One can
 % identify this group as the group of all $\Upsilon$ of the form
 % from Proposition 2 and having integer entries and determinant one.
 % \end{remark}
  Now $dim_C H^1(W,D_W) = 12 +11m +6m^2$ and $H^0(W,D_W)$ is
  generated by the elements in the center of $X$. To count the
independent parameters in the deformation space we need the
following:

 \begin{lemma}  All elements in the basis of $H^1(W,D_W)$ from
  (43) are invariant under hypercomplex transformations.
\end{lemma}
  \bproof  The elements in $H^0(W,D_W)$ are described in Lemma
  5. In particular the whole space is spanned by linear combinations
  of the vectors in the center $(Z,E_1,E_2,E_3)$ with coefficients
  depending on the fiber coordinate $\mu$. More precisely
  $H^0(W,D_W) = span\{W_k^{m+1} = I_kZ-iI_{\vec{a} }I_kZ\}$ at the
  point $\vec{a}$. Then $W_k^{m+1}$ are expressed locally as in (36)
  in terms of the local vector fields $\partial_i^{m+1}$. The
  $\partial_i^{m+1}$ themselves are defined in (34). In particular
  there is no $\frac{d}{d\mu}$ involved. The space $H^1(W,D_W)$ is a
  subspace of the span of
  $\lambda_1^k\lambda_2^{2-k}V_i^a\overline{\Omega}_j^b$. Again as
  in the calculation for the Nijenhuis bracket the elements are
  expressed locally in terms of
  $\frac{\mu^k}{|\mu|^2+1}\partial_i^{\alpha}\overline{\sigma}_j^{\beta}$
  for $k=0,1,2$. Here $\mu$ is inessential too.

   Now we have to check that ${\cal L}_X(A) = 0$ for $X \in
  H^0$ and $A \in H^1$. The formula we need is:
  $$
  {\cal L}_X (\alpha \otimes Y) = ({\cal L}_X\alpha)\otimes Y +
  \alpha \otimes ({\cal L}_X Y)
  $$

   Now we take $X$ to be any vector which is linear combination of
  $(Z,E_1,E_2,E_3)$ with coefficients depending on $\mu$. Then
  ${\cal L}_X V_i^a = 0$ for any $V_i^a$ - here $\mu$ is a constant
  for the differentiations. Moreover using

  $$
  {\cal L}_{\partial_i^{m+1} }
  (\overline{\sigma}_j^{\alpha})(\overline{\partial}_k^{\beta}) =
  \partial_i^{m+1}(\overline{\sigma}_j^{\alpha}(\overline{\partial}_k^{\beta})) -
  \overline{\sigma}_j^{\alpha}([ \partial_i^{m+1},
  \overline{\partial}_k^{\beta}]) = 0
  $$
  we see that
   ${\cal L}_X \overline{\Omega}_j^b(Y) = 0$ for $X$ in $H^0$
   because $Y$ is a combination of vectors $\partial_k^{\beta}$.

\eproof

  So by the Lemma all small deformations of the structures are again
invariant structures by Cathelineau's theory and there are no
equivalent deformations in
  $H^1$.
  \begin{corollary}
   There are deformations of the hypercomplex structure on $X$ which
  do not arise from a deformation of the lattice $\Gamma$.
 \end{corollary}

This combined with Proposition \ref{proposition2} proves Theorem
\ref{theorem5} from the Introduction.

\begin{remark}
 We note a similar phenomena for the complex deformations of the 3-dimensional complex Heisenberg
group $G$. The local moduli space around its canonical complex
structure is calculated by \cite{Nak} and is 6 dimensional. The
connected component of the invariant complex structures at this
point is again 6-dimensional, but the orbit of the complex
structure under the action of the automorphism group of $G$ is
only 2 dimensional \cite{Sal}. Since every two lattices are
equivalent under an automorphism of $G$, it follows that there are
deformations of the complex structure, which do not arise from
deformations of the lattice, unlike the case of the complex torus.
This also led Nakamura \cite{Nak} to the conclusion that a small
deformation of a complex parallelizable manifold is not necessary
complex parallelizable. What we proved in this section is that
similar fenomena holds for the small deformations of the
hypercomplex structure of $X$.
\end{remark}

\begin{remark} One could notice that above we included identifications of
the structures which arise from affine transformations. This
identifies additionally some structures in the orbit of the
automorphism group, since it acts only by linear transformations.
So the number of effective parameters for the small deformations
is less then the dimension of the orbit. Similar fenomena appears
in the complex deformations of $H_{2n+1}\times R/\Gamma$ as
described in \cite{GPP}. Due to the translation factor, there are
no "off-diagonal" deformations. The difference in the hypercomplex
case is that this amounts to all identifications because of the
properties of hypercomplex automorphisms of ${\bf H}^n$. In the
complex deformations case, this follows from the Kuranishi theory
\cite{MPPS}.

\end{remark}

        \end{document}